\documentclass[11pt,a4paper]{article}

\usepackage{float}
\usepackage[T1]{fontenc}
\usepackage[margin=2.7cm]{geometry}
\usepackage{amsmath,amssymb,amsthm,mathtools,bm}
\usepackage{graphicx,booktabs,array,enumitem}
\usepackage{tikz}
\usetikzlibrary{arrows.meta,positioning,fit,backgrounds,calc}
\usepackage[hidelinks]{hyperref}
\theoremstyle{plain}
\newtheorem{theorem}{Theorem}
\newtheorem{proposition}{Proposition}
\newtheorem{lemma}{Lemma}
\newtheorem{corollary}{Corollary}
\theoremstyle{definition}

\theoremstyle{remark}
\newtheorem{remark}{Remark}
\newtheorem{example}{Example}

\numberwithin{equation}{section}
\newcommand{\M}{\mathcal M}
\newcommand{\R}{\mathbb R}
\newcommand{\T}{\mathrm T}
\DeclareMathOperator*{\argmin}{arg\,min}

\title{Moving Cut Loci and Gaussian Fluctuations of Intrinsic Huber Means on Compact Riemannian Manifolds}

\author{
  Rui Qiu \qquad Yiwei Yao \qquad Zhou Yu
  \\[0.6em]
  {\small School of Statistics, East China Normal University}
}

\date{}
\begin{document}

\maketitle

\begin{abstract}
Moving cut loci can change the limiting covariance of sample Huber means even when the cut locus of the population mean has probability zero.
On a closed connected Riemannian manifold, we prove that the averaged logarithm field is $C^1$ near a point $m$ whenever the sampling law has a continuous density near its cut locus.
The result holds for every continuous weight depending only on the observation, and the covariant derivative includes an explicit integral over the regular cut locus.
When $m$ is a Huber mean for a fixed positive threshold, assume that the distance to $m$ has no atom at that threshold.
The Hessian of the population risk then equals the expected pointwise Hessian minus a positive semidefinite cut locus correction.
If $m$ is unique and this Hessian is positive definite, every measurable selection of global minimizers for i.i.d.\ samples satisfies a tangent-space central limit theorem.
The cut loci of points in a ball of radius $\delta$ sweep out a region of volume $O(\delta)$, uniformly in the centre.
These results establish the volume and integrability conditions previously assumed in a central limit theorem for intrinsic Fr\'echet means.
At each fixed centre, the volume measure on this region, divided by $\delta$, converges weakly to an explicit measure on the regular cut locus.
The Huber correction vanishes exactly when the region swept out about $m$ has probability $o(\delta)$.
\end{abstract}

\section{Introduction}
\label{sec:introduction}

For manifold-valued observations, the Fr\'echet mean extends
the Euclidean mean by minimizing the expected squared geodesic
distance \cite{Frechet1948}.
Its geometric foundations and computational methods have been
studied extensively \cite{Karcher1977,Afsari2011,Pennec2006},
alongside consistency and central limit theory
\cite{BhattacharyaPatrangenaru2003,BhattacharyaPatrangenaru2005,
KendallLe2011,BhattacharyaLin2017}.
The squared-distance criterion assigns larger score magnitudes
to observations farther from the centre, motivating robust
alternatives such as the Riemannian geometric median
\cite{FletcherEtAl2009,Yang2010}.
Huber's loss \cite{Huber1964} is quadratic near the centre
and linear beyond a prescribed threshold, thereby capping
the radial score.
A Huber mean minimizes the expected Huber loss of the geodesic
distance, providing a robust alternative within the Fr\'echet
framework.

Lee and Jung \cite{LeeJung2026} developed the statistical theory
of Huber means on Riemannian manifolds, establishing existence,
conditions for uniqueness, strong consistency of sample mean sets,
and asymptotic normality under suitable regularity assumptions.
Their central limit theorem expresses the limiting covariance in
terms of the expected pointwise Hessian.
When the support of the distribution meets the cut
locus of the population Huber mean, does this expected Hessian
still equal the Hessian of the population risk?

Cut locus stability holds on compact Riemannian manifolds
\cite{EltznerEtAl2021}, but moving the base point can still
produce a first-order contribution to the mean score even
when a fixed cut locus has zero probability.
On the unit circle, moving the base point through
an angle $t$ moves its antipode through the same angle. For each
observation in the interior of the swept arc, the minimizing
geodesic switches from one direction around the circle to the other.
If the density is continuous and positive at the original antipode,
the probability of this arc is of order $|t|$ as $t\to0$
\cite{HotzHuckemann2015}.
Hotz, Le and Wood \cite{HotzLeWood2024} obtained a central
limit theorem for intrinsic Fr\'echet means with an explicit
cut locus correction, assuming a continuous density near
the cut locus, a linear swept-volume bound and integrability
of the derivative of the logarithm map.
Aveni \cite[Chapter~4]{Aveni2025Thesis} treated a class of
twice continuously differentiable radial losses under
related volume and integrability assumptions.
This leaves two questions: whether these additional
regularity conditions follow from the geometric and local
density assumptions, and how to treat Huber's loss, which
is continuously differentiable but not twice differentiable
at the threshold.

Our main analytic contribution is a differentiation theorem
for weighted logarithm fields on smooth compact connected
Riemannian manifolds without boundary.
Proposition~\ref{prop:field-difference-hotz} shows that a
continuous volume density near the cut locus of a base point
suffices for local $C^1$ regularity of the averaged logarithm
field for every fixed continuous weight depending only on
the observation.
Its covariant derivative consists of regular and cut locus
integrals, both absolutely convergent and continuous in the
base point.
Theorem~\ref{thm:sweep-volume} gives the uniform
linear swept-volume bound in condition~\textup{(iv)} of
\cite[Theorem~1]{HotzLeWood2024} from the geometry alone.
Under the local density assumption, Proposition~\ref{prop:field-difference-hotz}
with unit weight gives the integrability of the derivative of the
logarithm map required in condition~\textup{(v)}.
To apply the proposition to the Huber score, we first hold the
radial weight fixed at its value at the population Huber mean,
then expand the change in the weight as the base point moves.

For a fixed threshold $c>0$, let $F_c$ be the population
Huber risk and $m_c$ a population minimizer.
Under the same local density assumption and zero probability
on the threshold sphere when $c<\operatorname{diam}(\M)$,
Theorem~\ref{thm:huber-score-linearization} and
Corollary~\ref{cor:huber-risk-quadratic} identify the
population Hessian as
$$
\Lambda_{\mu,c}
=\Lambda_{\mu,c}^{\mathrm{reg}}-J_{\mu,c}(m_c),
$$
where $\Lambda_{\mu,c}^{\mathrm{reg}}$ is the expected
pointwise Hessian and $J_{\mu,c}(m_c)$ is an explicit
positive semidefinite cut locus correction.
If $m_c$ is unique and $\Lambda_{\mu,c}$ is strictly positive
definite, Theorem~\ref{thm:huber-clt} gives a $\sqrt n$
central limit theorem for every measurable selection of
sample global minimizers from independent and identically
distributed observations, with limiting covariance
$\Lambda_{\mu,c}^{-1}\Sigma_c\Lambda_{\mu,c}^{-1}$,
where $\Sigma_c$ is the covariance of the individual Huber score
at $m_c$.

The swept-volume theorem also yields a geometric criterion
for the cut locus correction to vanish.
Let $C_x$ denote the cut locus of $x$ and write
$\mathcal A_\delta(m)=\bigcup_{x\in B_\delta(m)}C_x$.
On every closed connected Riemannian manifold,
Theorem~\ref{thm:sweep-volume} gives
$\operatorname{vol}(\mathcal A_\delta(m))=O(\delta)$ uniformly in $m$.
The linear order is sharp in general.
Banerjee, Banerjee and Das
\cite[Lemma~4]{BanerjeeBanerjeeDas2026} obtained the corresponding
volume bound on real Grassmannians.
For each fixed $m$, the rescaled volume measure
$\delta^{-1}\mathbf{1}_{\mathcal A_\delta(m)}\,\mathrm d\operatorname{vol}$
converges weakly to an explicit measure on the regular cut locus of $m$.
The total mass of this limit gives the exact first-order coefficient of the
swept volume. Under the local continuous-density assumption, the
corresponding probability expansion gives
$$
J_{\mu,c}(m_c)=0\Longleftrightarrow
\mu\bigl(\mathcal A_\delta(m_c)\bigr)=o(\delta)
 \text{ as }\delta\downarrow0
$$
(Corollary~\ref{cor:sweep-correction-vanishing}). Under the threshold condition stated above, this equivalence characterizes precisely when the expected pointwise Hessian equals the population Hessian: the probability swept out by nearby cut loci must be negligible to first order.
Explicit examples on the circle, the flat two-dimensional
torus and $\mathrm{SO}(3)$ show that omitting the correction
can change the limiting covariance.

\section{Main Results}
\label{sec:main-results}
\subsection{Setup and notation}
\label{subsec:notation}

Let $(\M,g)$ be a smooth, compact, connected Riemannian manifold 
of dimension $d\ge1$ without boundary, and let $\mu$ be a Borel
probability measure on $\M$. We write $d(\cdot,\cdot)$ for the
Riemannian distance, $D=\operatorname{diam}(\M)$ for the diameter,
and $\operatorname{vol}$ for the Riemannian volume measure.
The metric induces an inner product $\langle\cdot,\cdot\rangle$
and a norm $\|\cdot\|$ on each tangent space, and
$\|\cdot\|_{\mathrm{op}}$ denotes the induced operator norm.

For $x\in\M$, let $\T_x\M$ be the tangent space at $x$,
$C_x=\operatorname{Cut}(x)$ the cut locus of $x$, and
$\operatorname{Exp}_x:\T_x\M\to\M$ the exponential map.
Let $\operatorname{Log}_x:\M\setminus C_x\to\T_x\M$
denote the Riemannian logarithm map, namely the inverse of
$\operatorname{Exp}_x$ restricted to its injectivity domain.
For sufficiently close $x,z\in\M$, we denote parallel transport
along their unique minimizing geodesic by
$$
\Pi_{x,z}:\T_x\M\longrightarrow\T_z\M.
$$

Let $D_V$ denote the Levi--Civita covariant derivative in the
direction $V$. For a $C^1$ vector field $Z$, we regard $DZ(x)$
as the linear operator on $\T_x\M$ given by
$DZ(x)V=D_VZ(x)$. For vector fields depending on two points,
such as $g_c(x,y)$ below, differentiation always acts on the
first variable with the second held fixed. For $n\in\T_x\M$, we define the operator $n\otimes n$ by
$$
(n\otimes n)V=\langle n,V\rangle n,
\qquad V\in\T_x\M.
$$

Fix a cutoff parameter $c\in(0,\infty)$. The Huber loss and its
associated radial weight are
$$
\rho_c(r)=
\begin{cases}
r^2, & 0\le r\le c,\\
2cr-c^2, & r>c,
\end{cases}
\qquad
a_c(r)=
\begin{cases}
2, & 0\le r\le c,\\[1mm]
\dfrac{2c}{r}, & r>c.
\end{cases}
$$
For $y\notin C_x$, define the Huber score and the covariant
derivative of the logarithm map by
$$
g_c(x,y)=a_c(d(x,y))\operatorname{Log}_x(y),\qquad
H(x\mid y)V=D_V\operatorname{Log}_x(y),\qquad V\in\T_x\M.
$$
For $y\in C_x$, we set
$\operatorname{Log}_x(y)=g_c(x,y)=0$ and $H(x\mid y)=0$.
At $y=x$, we have
$H(x\mid x)=-\operatorname{Id}_{\T_x\M}$.

The population Huber risk and the mean score field are
$$
F_c(x)=\int_{\M}\rho_c(d(x,y))\,\mathrm d\mu(y),\qquad
G_c(x)=\int_{\M\setminus C_x}g_c(x,y)\,\mathrm d\mu(y).
$$
A population Huber mean is any global minimizer of $F_c$.
The score satisfies $\|g_c(x,y)\|\le2c$ for all $x,y\in\M$.

We next introduce the geometric quantities appearing in the
cut locus correction. Let $H_x$ denote the set of points
$y\in C_x$ joined to $x$ by exactly two minimizing geodesics,
such that $y$ is not conjugate to $x$ along either geodesic.
Set
$$
Q_x=C_x\setminus H_x.
$$
We refer to $H_x$ and $Q_x$ as the regular and exceptional
parts of the cut locus, respectively; see
\cite[Lemma~1 and the subsequent discussion]{HotzLeWood2024}.

Let $\mathcal H^{d-1}$ denote the $(d-1)$-dimensional Hausdorff
measure induced by the Riemannian distance.
By \cite[Theorem~2.7]{MMU2014}, $\mathcal H^{d-1}(Q_x)=0$ for $d\ge2$.
When $d=1$, the manifold is isometric to a circle and $Q_x=\varnothing$.

For each $y\in H_x$, there are a neighborhood $V_y$ of $y$
and two smooth local distance branches
$\varphi_{1,y}(\,\cdot\mid x)$ and
$\varphi_{2,y}(\,\cdot\mid x)$ such that
$$
d(x,z)=\min\{
\varphi_{1,y}(z\mid x),
\varphi_{2,y}(z\mid x)
\},
\qquad z\in V_y.
$$
Their equality set is $V_y\cap H_x$, and the gradient of
their difference is nonzero at $y$; see
\cite[Lemma~2 and Corollary~2]{HotzLeWood2024}. Define
$$
\begin{aligned}
\kappa(y\mid x)
&=
\left\|
\left.\operatorname{grad}_z\bigl(
\varphi_{1,y}(z\mid x)-\varphi_{2,y}(z\mid x)
\bigr)\right|_{z=y}
\right\|,\\
n(y\mid x)
&=
\frac{
\left.\operatorname{grad}_z\bigl(
\varphi_{1,y}(z\mid x)-\varphi_{2,y}(z\mid x)
\bigr)\right|_{z=y}
}{\kappa(y\mid x)}.
\end{aligned}
$$
Thus $\kappa(y\mid x)>0$, and $n(y\mid x)\in\T_y\M$ is
a unit normal to $H_x$ at $y$.
Finally, the reciprocity relations
$y\in C_x\Longleftrightarrow x\in C_y,
y\in H_x\Longleftrightarrow x\in H_y$
hold; see \cite{HotzLeWood2024}. In particular, for $y\in H_x$,
the quantities $\kappa(x\mid y)$ and $n(x\mid y)\in\T_x\M$
are defined by interchanging the roles of $x$ and $y$.
\subsection{Regularity and linearization of weighted logarithm vector fields}
\begin{proposition}
\label{prop:field-difference-hotz}
Fix $m\in\M$. Suppose that there exists an open neighborhood
$O$ of $C_m$ such that $\mu|_O$ is absolutely continuous with
respect to $\operatorname{vol}|_O$ and the corresponding
Radon--Nikodym derivative admits a continuous version
$\psi:O\to[0,\infty)$. Thus,
$$
\mu|_O=\psi\,\operatorname{vol}|_O.
$$
Let $w:\M\to\R$ be a fixed continuous function and define
$$
L_{\mu,w}(x)
=\int_{\M\setminus C_x}w(y)\operatorname{Log}_x(y)\,\mathrm d\mu(y).
$$
Then there exists a normal neighborhood $U$ of $m$ such that
the following statements hold.

\textup{(i)} For every $x\in U$, we have $C_x\subset O$ and
$\mu(C_x)=0$.

\textup{(ii)} For $x\in U$, define
\begin{equation}
\label{eq:weighted-A-definition}
A_{\mu,w}(x)
=\int_{\M\setminus C_x}w(y)H(x\mid y)\,\mathrm d\mu(y)
\end{equation}
and
\begin{equation}
\label{eq:weighted-J-definition}
J_{\mu,w}(x)
=\int_{H_x}w(y)d(x,y)
\frac{\kappa(x\mid y)^2}{\kappa(y\mid x)}\psi(y)
\bigl(n(x\mid y)\otimes n(x\mid y)\bigr)
\,\mathrm d\mathcal H^{d-1}(y).
\end{equation}
Both operator-valued integrals converge absolutely, and
$A_{\mu,w}$ and $J_{\mu,w}$ are continuous on $U$.
Moreover, $L_{\mu,w}\in C^1(U;\T\M)$ with
\begin{equation}
\label{eq:weighted-continuous-derivative}
DL_{\mu,w}(x)=A_{\mu,w}(x)+J_{\mu,w}(x).
\end{equation}

\textup{(iii)} As $x\to m$,
\begin{equation}
\label{eq:weighted-full-linearization}
\Pi_{x,m}L_{\mu,w}(x)-L_{\mu,w}(m)=\bigl[A_{\mu,w}(m)+J_{\mu,w}(m)\bigr]
\operatorname{Log}_m(x)+o(d(x,m)).
\end{equation}
\end{proposition}

The proof is given in Section~\ref{sec:proof-weighted-linearization}.
\subsection{Linearization of the Huber score vector field}

Fix a base point $m_c\in\M$ and write
$$
S_c(m_c)=\{y\in\M:d(m_c,y)=c\}.
$$
Whenever the following integral converges absolutely, define
the operator $\Lambda_{\mu,c}^{\mathrm{reg}}$ on $\T_{m_c}\M$ by
$$
\Lambda_{\mu,c}^{\mathrm{reg}}V
=-\int_{\M\setminus(C_{m_c}\cup S_c(m_c))}
D_Vg_c(m_c,y)\,\mathrm d\mu(y),
\qquad V\in\T_{m_c}\M.
$$
If $\mu$ has a continuous density $\psi$ on a neighborhood of
$C_{m_c}$, apply Proposition~\ref{prop:field-difference-hotz}
with base point $m_c$ and fixed weight $w(y)=a_c(d(m_c,y))$,
and define
$$
\begin{aligned}
J_{\mu,c}(m_c)
&:=J_{\mu,w}(m_c)\\
&=\int_{H_{m_c}}2\min\{d(m_c,y),c\}
\frac{\kappa(m_c\mid y)^2}{\kappa(y\mid m_c)}\psi(y)
\bigl(n(m_c\mid y)\otimes n(m_c\mid y)\bigr)
\,\mathrm d\mathcal H^{d-1}(y).
\end{aligned}
$$
Set
$$
\Lambda_{\mu,c}
=\Lambda_{\mu,c}^{\mathrm{reg}}-J_{\mu,c}(m_c).
$$

\begin{theorem}
\label{thm:huber-score-linearization}
Fix $c\in(0,\infty)$. Suppose that the following conditions hold.

\textup{(i)} $\M$ is a smooth, compact, connected
$d$-dimensional Riemannian manifold without boundary.

\textup{(ii)} The point $m_c$ is a local minimizer of $F_c$.

\textup{(iii)} On some neighborhood of $C_{m_c}$, the measure
$\mu$ is absolutely continuous with respect to Riemannian
volume, and its Radon--Nikodym derivative admits a continuous
version $\psi$.

\textup{(iv)} If $0<c<D$, then $\mu(S_c(m_c))=0$.

Then $\Lambda_{\mu,c}^{\mathrm{reg}}$ and $J_{\mu,c}(m_c)$
are well defined, with absolutely convergent defining integrals.
Moreover, $\Lambda_{\mu,c}^{\mathrm{reg}}$, $J_{\mu,c}(m_c)$,
and $\Lambda_{\mu,c}$ are self-adjoint, and
$$
\mu(C_{m_c}\cup S_c(m_c))=0,
\qquad G_c(m_c)=0.
$$
As $x\to m_c$,
$$
\Pi_{x,m_c}G_c(x)
=-\Lambda_{\mu,c}\operatorname{Log}_{m_c}(x)
+o(d(x,m_c)).
$$
If $H_{m_c}=\varnothing$, then $J_{\mu,c}(m_c)=0$, so
$\Lambda_{\mu,c}=\Lambda_{\mu,c}^{\mathrm{reg}}$.
\end{theorem}

A self-adjoint operator $A$ is called strictly positive definite
if $\langle V,AV\rangle>0$ for every nonzero $V$.
We use $\operatorname{Hess}^{g}f$ for the Riemannian Hessian
as a bilinear form:
$$
\operatorname{Hess}^{g}f(V,W)
=\langle D_V\operatorname{grad}f,W\rangle.
$$

\begin{corollary}
\label{cor:huber-risk-quadratic}
Under the assumptions of Theorem~\ref{thm:huber-score-linearization},
$F_c$ belongs to $C^1(U)$ for some neighborhood $U$ of $m_c$,
and its gradient is Fr\'echet differentiable at $m_c$.
In particular, its Riemannian Hessian at $m_c$ exists and satisfies
$$
\operatorname{Hess}^{g}F_c(m_c)(V,W)
=\langle\Lambda_{\mu,c}V,W\rangle,
\qquad V,W\in\T_{m_c}\M.
$$
As $u\to0$ in $\T_{m_c}\M$,
$$
F_c(\operatorname{Exp}_{m_c}u)-F_c(m_c)
=\frac12\langle u,\Lambda_{\mu,c}u\rangle
+o(\|u\|^2).
$$
If $\Lambda_{\mu,c}$ is strictly positive definite, then there
exist $a>0$ and $\delta>0$ such that
$$
F_c(x)-F_c(m_c)\ge a\,d(x,m_c)^2
$$
whenever $d(x,m_c)<\delta$.
\end{corollary}

The proofs of Theorem~\ref{thm:huber-score-linearization} and
Corollary~\ref{cor:huber-risk-quadratic} are given in
Section~\ref{sec:proof-linearization}.

\subsection{Central limit theorem}
\label{subsec:huber-clt}

Let $X_1,X_2,\ldots$ be independent and identically distributed
with law $\mu$, and define
$$
F_{n,c}(x)=\frac1n\sum_{i=1}^n\rho_c(d(x,X_i)).
$$
The set of sample minimizers is
$$
E_{n,c}=\argmin_{x\in\M}F_{n,c}(x).
$$
Since $\M$ is compact and $F_{n,c}$ is continuous, $E_{n,c}$
is nonempty and compact.
The sample objective is jointly continuous in
$(x,X_1,\ldots,X_n)$, so the measurable maximum theorem
ensures that $E_{n,c}$ admits a measurable selection.
Let $\widehat m_{n,c}$ be any such selection.

For a random vector $Y\in\T_{m_c}\M$ with finite second moment,
define its covariance operator
$\operatorname{Cov}(Y):\T_{m_c}\M\to\T_{m_c}\M$ by
$$
\operatorname{Cov}(Y)V
=\mathbb E\left[
\langle Y-\mathbb E Y,V\rangle(Y-\mathbb E Y)
\right],
\qquad V\in\T_{m_c}\M.
$$
For a self-adjoint positive semidefinite operator $\Gamma$
on $\T_{m_c}\M$, let $N_{\T_{m_c}\M}(0,\Gamma)$ denote the
Gaussian distribution with mean zero and covariance operator
$\Gamma$. Set
$$
\Sigma_c=\operatorname{Cov}\{g_c(m_c,X_1)\}.
$$
This covariance is well defined because
$\|g_c(m_c,X_1)\|\le2c$.

Fix a normal neighborhood $U_0=\operatorname{Exp}_{m_c}(B_{r_0}(0)), 0<r_0<\operatorname{inj}(m_c)$, where $B_{r_0}(0)$ is the open ball of radius $r_0$ centered at the origin in $\T_{m_c}\M$, and $\operatorname{inj}(m_c)$
is the injectivity radius at $m_c$. Define the zero extension of
$\operatorname{Log}_{m_c}|_{U_0}$ to $\M$ by
$$
\widetilde{\operatorname{Log}}_{m_c}(x)
=\begin{cases}
\operatorname{Log}_{m_c}(x),&x\in U_0,\\
0,&x\notin U_0.
\end{cases}
$$

\begin{theorem}
\label{thm:huber-clt}
Suppose that the assumptions of
Theorem~\ref{thm:huber-score-linearization} hold.
Assume further that $m_c$ is the unique global minimizer of
$F_c$ and that $\Lambda_{\mu,c}$ is strictly positive definite.
Then every measurable selection $\widehat m_{n,c}$ satisfies
$$
\sqrt n\,\widetilde{\operatorname{Log}}_{m_c}(\widehat m_{n,c})
\overset{d}{\longrightarrow}
N_{\T_{m_c}\M}\left(
0,\Lambda_{\mu,c}^{-1}\Sigma_c\Lambda_{\mu,c}^{-1}
\right).
$$
\end{theorem}

The proof is given in Section~\ref{sec:proof-clt}.

\subsection{Volume asymptotics for moving cut loci and vanishing of the correction}
\label{subsec:sweep-volume-results}

For $m\in\M$ and $\delta>0$, define
$$
\mathcal A_\delta(m)=\bigcup_{x\in B_\delta(m)}C_x,
$$
where $B_\delta(m)=\{x\in\M:d(x,m)<\delta\}$.

\begin{theorem}
\label{thm:sweep-volume}
Let $(\M,g)$ be a smooth, compact, connected Riemannian
manifold of dimension $d\ge1$ without boundary.
Then the following statements hold.

\textup{(i)} There exist constants $C>0$ and $\delta_0>0$,
depending only on $(\M,g)$, such that
$$
\operatorname{vol}(\mathcal A_\delta(m))\le C\delta,
\qquad m\in\M,\quad 0<\delta<\delta_0.
$$

\textup{(ii)} For every fixed $m\in\M$ and every $f\in C(\M)$,
$$
\lim_{\delta\downarrow0}\frac1\delta
\int_{\mathcal A_\delta(m)}f(y)\,\mathrm d\operatorname{vol}(y)
=2\int_{H_m}
f(y)\frac{\kappa(m\mid y)}{\kappa(y\mid m)}
\,\mathrm d\mathcal H^{d-1}(y),
$$
and the integral on the right converges absolutely.
In particular,
$$
\operatorname{vol}(\mathcal A_\delta(m))
=2\delta\int_{H_m}
\frac{\kappa(m\mid y)}{\kappa(y\mid m)}
\,\mathrm d\mathcal H^{d-1}(y)+o(\delta).
$$
\end{theorem}

The proof of Theorem~\ref{thm:sweep-volume} is given in
Appendix~\ref{app:sweep-volume-proof}.

\begin{corollary}
\label{cor:sweep-probability}
Let $\mu$ be a Borel probability measure on $\M$ and fix $m\in\M$.
Suppose that there exists an open neighborhood $O$ of $C_m$
such that
$$
\mu|_O=\psi\,\operatorname{vol}|_O,
\qquad \psi\in C(O),\quad \psi\ge0.
$$
Then, for every $w\in C(\M)$, as $\delta\downarrow0$,
$$
\int_{\mathcal A_\delta(m)}w(y)\,\mathrm d\mu(y)
=2\delta\int_{H_m}
w(y)\psi(y)\frac{\kappa(m\mid y)}{\kappa(y\mid m)}
\,\mathrm d\mathcal H^{d-1}(y)+o(\delta).
$$
In particular,
$$
\mu(\mathcal A_\delta(m))
=2\delta\int_{H_m}
\psi(y)\frac{\kappa(m\mid y)}{\kappa(y\mid m)}
\,\mathrm d\mathcal H^{d-1}(y)+o(\delta)=O(\delta).
$$
\end{corollary}

\begin{proof}
By Lemma~\ref{lem:cut-time-closedness}\textup{(iii)},
$\mathcal A_\delta(m)$ is contained in any prescribed open
neighborhood of $C_m$ for all sufficiently small $\delta$.
Choose $\eta\in C^\infty(\M)$ with
$\operatorname{supp}\eta\Subset O$ such that $\eta=1$ on a
neighborhood of $C_m$.
Extend $\eta w\psi$ by zero outside $O$ to obtain a function
in $C(\M)$.
For all sufficiently small $\delta$,
$$
\int_{\mathcal A_\delta(m)}w\,\mathrm d\mu
=\int_{\mathcal A_\delta(m)}
\eta w\psi\,\mathrm d\operatorname{vol}.
$$
The result follows from Theorem~\ref{thm:sweep-volume}\textup{(ii)}.
\end{proof}
\begin{corollary}
\label{cor:sweep-correction-vanishing}
Suppose that the assumptions of
Corollary~\ref{cor:sweep-probability} hold with $m=m_c$.
Then
$$
J_{\mu,c}(m_c)=0\Longleftrightarrow
\mu(\mathcal A_\delta(m_c))=o(\delta)\ \text{as }\delta\downarrow0.
$$
\end{corollary}

\begin{proof}
By Corollary~\ref{cor:sweep-probability},
$\mu(\mathcal A_\delta(m_c))=o(\delta)$ as $\delta\downarrow0$
if and only if
$$
\int_{H_{m_c}}
\psi(y)\frac{\kappa(m_c\mid y)}{\kappa(y\mid m_c)}
\,\mathrm d\mathcal H^{d-1}(y)=0.
$$
Since $\psi\ge0$ and the ratio is strictly positive on
$H_{m_c}$, this holds if and only if $\psi=0$ almost
everywhere on $H_{m_c}$ with respect to $\mathcal H^{d-1}$.

The operator $J_{\mu,c}(m_c)$ is positive semidefinite and
$$
\operatorname{tr}J_{\mu,c}(m_c)
=\int_{H_{m_c}}
2\min\{d(m_c,y),c\}
\frac{\kappa(m_c\mid y)^2}{\kappa(y\mid m_c)}
\psi(y)\,\mathrm d\mathcal H^{d-1}(y).
$$
The coefficient multiplying $\psi$ is strictly positive
on $H_{m_c}$. Thus $\operatorname{tr}J_{\mu,c}(m_c)=0$ if and only if
$\psi=0$ almost everywhere on $H_{m_c}$ with respect to
$\mathcal H^{d-1}$. A positive semidefinite operator is zero if and only if
its trace is zero, which proves the equivalence.
\end{proof}

\begin{remark}
\label{rem:sweep-sharpness}
On the flat torus $\M=(\R/2\pi\mathbb Z)^d$ with the product metric,
for $0<\delta<\pi$, the set $\mathcal A_\delta(m)$ is exactly
the union of the strips of width $2\delta$ around the $d$
coordinate cut hypersurfaces. Hence
$$
\operatorname{vol}(\mathcal A_\delta(m))
=(2\pi)^d-(2\pi-2\delta)^d
=2d(2\pi)^{d-1}\delta+O(\delta^2).
$$
Thus, in general, the linear upper bound in
Theorem~\ref{thm:sweep-volume} cannot be replaced by
$o(\delta)$ or by $O(\delta^{1+\varepsilon})$ for any
$\varepsilon>0$.
\end{remark}

\begin{remark}
\label{rem:weighted-assumptions}
When $c\ge D$, the Huber loss reduces to squared distance
and the limiting covariance agrees with that of
\cite[Theorem~2]{HotzLeWood2024} whenever both theorems apply.
Condition~\textup{(iv)} of Theorem~1 in
\cite{HotzLeWood2024} follows from
Theorem~\ref{thm:sweep-volume}\textup{(i)}.
Under our local density assumption,
Proposition~\ref{prop:field-difference-hotz} with $w\equiv1$
also gives the integrability required in condition~\textup{(v)}
of that theorem.

In normal coordinates induced by an orthonormal basis of
$\T_{m_c}\M$, the expected pointwise Hessian matrix $H_c$
in \cite[Theorem~4]{LeeJung2026} represents
$\Lambda_{\mu,c}^{\mathrm{reg}}$.
By Corollary~\ref{cor:huber-risk-quadratic}, the matrix of
the population Hessian in the same basis is that of
$\Lambda_{\mu,c}^{\mathrm{reg}}-J_{\mu,c}(m_c)$.
The correction $J_{\mu,c}(m_c)$ can be nonzero even when
$\mu(C_{m_c})=0$. When the correction vanishes and both
theorems apply, the limiting covariance formulas agree.
Assumption~\textup{(A3)} in \cite{LeeJung2026} requires
almost sure uniqueness of the sample Huber mean for every
sample size, whereas Theorem~\ref{thm:huber-clt} allows
any measurable selection from the sample minimizer set.
\end{remark}

\section{Examples}
\label{sec:example}

\subsection{A smooth density on \texorpdfstring{$S^1$}{S1}}

Let $\M=S^1$ be the unit circle with arc-length coordinate
$\theta\in(-\pi,\pi]$. We represent tangent vectors and
linear operators by their scalar components in the unit
coordinate frame. Fix $0<c<\pi$ and consider the density
$$
\psi_a(\theta)=\frac{1+a\cos\theta}{2\pi},
\qquad 0<a<1.
$$
Write $\mu_a$ for the corresponding probability measure
and let $\Theta$ denote the angular coordinate of an
observation with law $\mu_a$.
The density is positive at the cut point $\pi$ of $0$.

\begin{example}
\label{ex:circle-explicit}
For the distribution $\mu_a$, the point $m_c=0$ is the
unique global minimizer of $F_c$ and
$$
\Lambda_{\mu_a,c}
=\frac{2a(c+\sin c)}{\pi}>0,\qquad
\Sigma_c
=\frac4\pi
\left[\pi c^2-\frac23c^3
+2a(c\cos c-\sin c)\right].
$$
Consequently,
$$
\sqrt n\,\widetilde{\operatorname{Log}}_0(\widehat m_{n,c})
\overset{d}{\longrightarrow}
N\left(
0,
\frac{\pi\left[3\pi c^2-2c^3
+6a(c\cos c-\sin c)\right]}
{3a^2(c+\sin c)^2}
\right).
$$
For $0<\delta<\pi$, we also have $\operatorname{vol}(\mathcal A_\delta(0))=2\delta$.
Thus the first-order volume coefficient in Theorem~\ref{thm:sweep-volume} is $2$.
\end{example}

\begin{proof}
Extend $q_c(t)=\rho_c(|t|)$ from $[-\pi,\pi]$ to a
$2\pi$-periodic function. Its evenness and the identity
$$
\int_0^\pi q_c(t)\cos t\,\mathrm dt=-2(c+\sin c)
$$
give
$$
F_c(x)
=\frac1{2\pi}\int_{-\pi}^{\pi}
q_c(x-\theta)(1+a\cos\theta)\,\mathrm d\theta
=\pi c-c^2+\frac{c^3}{3\pi}
-\frac{2a(c+\sin c)}{\pi}\cos x.
$$
Hence $0$ is the unique global minimizer and
$$
\operatorname{Hess}^{g}F_c(0)(V,W)
=\frac{2a(c+\sin c)}{\pi}\langle V,W\rangle,
\qquad V,W\in\T_0\M.
$$

The pointwise derivative away from the cut locus and the
threshold points gives
$$
\Lambda_{\mu_a,c}^{\mathrm{reg}}
=2\mu_a(\{\theta:|\theta|<c\})
=\frac{2(c+a\sin c)}{\pi}.
$$
Since $H_0=C_0=\{\pi\}$ and
$\kappa(0\mid\pi)=\kappa(\pi\mid0)=2$, we obtain
$$
J_{\mu_a,c}(0)=4c\psi_a(\pi)
=\frac{2c(1-a)}{\pi}.
$$
Subtracting the correction from the regular term agrees
with the population Hessian computed above.

By symmetry the score has mean zero and
$ g_c(0,\Theta)=2\operatorname{sgn}(\Theta)\min\{|\Theta|,c\}$
almost surely. Therefore,
$$
\begin{aligned}
\Sigma_c
&=\frac4\pi\left[
\int_0^c t^2(1+a\cos t)\,\mathrm dt
+c^2\int_c^\pi(1+a\cos t)\,\mathrm dt
\right]\\
&=\frac4\pi\left[
\pi c^2-\frac23c^3+2a(c\cos c-\sin c)
\right].
\end{aligned}
$$
The measure $\mu_a$ has a smooth density and satisfies
$\mu_a(S_c(0))=0$.
Together with uniqueness and $\Lambda_{\mu_a,c}>0$,
these facts allow us to apply Theorem~\ref{thm:huber-clt}
and obtain the stated limiting distribution.

Finally, $C_x$ consists of the antipodal point of $x$.
Thus $\mathcal A_\delta(0)=B_\delta(\pi)$ has length
$2\delta$. Since $\mathcal H^0(H_0)=1$ and the ratio of
the two endpoint gradient norms is $1$,
Theorem~\ref{thm:sweep-volume} gives the same coefficient $2$.
\end{proof}

\begin{remark}
As $a\downarrow0$, the measures $\mu_a$ converge to the
uniform distribution. Both
$\Lambda_{\mu_a,c}^{\mathrm{reg}}$ and $J_{\mu_a,c}(0)$
converge to $2c/\pi$, so $\Lambda_{\mu_a,c}\to0$.
This agrees with the fact that the population risk is
constant under the uniform distribution and shows that
the cut locus correction can cancel the regular term.
\end{remark}

\subsection{The flat two-dimensional torus}

Let $\M=\mathbb T^2=\R^2/(2\pi\mathbb Z)^2$ carry the
standard flat metric. Represent its points by the fundamental
domain $(-\pi,\pi]^2$ and let $e_1,e_2$ be the standard
orthonormal coordinate frame. We represent operators by
matrices in this frame and write $\operatorname{Id}_2$
for the identity matrix.
Let $\mu_{\mathrm{Haar}}$ be the normalized Haar probability
measure on $\mathbb T^2$ and consider
$$
\mu_\varepsilon
=(1-\varepsilon)\mu_{\mathrm{Haar}}
+\varepsilon\delta_{(0,0)},
\qquad 0<\varepsilon<1,
$$
where $\delta_p$ denotes the Dirac measure at $p$.
Throughout this example we assume $0<c<\pi$.

\begin{example}
\label{ex:torus-explicit}
For the distribution $\mu_\varepsilon$, the point
$m_c=(0,0)$ is the unique global minimizer of $F_c$ and
$$
\Lambda_{\mu_\varepsilon,c}
=2\varepsilon\operatorname{Id}_2,\qquad
\Sigma_c
=2(1-\varepsilon)
\left(c^2-\frac{c^4}{8\pi}\right)\operatorname{Id}_2,
$$
and hence
$$
\sqrt n\,\widetilde{\operatorname{Log}}_{m_c}(\widehat m_{n,c})
\overset{d}{\longrightarrow}
N_{\T_{m_c}\M}\left(
0,
\frac{1-\varepsilon}{2\varepsilon^2}
\left(c^2-\frac{c^4}{8\pi}\right)\operatorname{Id}_2
\right).
$$
For $0<\delta<\pi$, we also have $\operatorname{vol}(\mathcal A_\delta(m_c))=8\pi\delta-4\delta^2$.
Thus the first-order volume coefficient in
Theorem~\ref{thm:sweep-volume} is $8\pi$.
\end{example}

\begin{proof}
Translation invariance and integration over the fundamental
domain give
$$
F_c(x)
=(1-\varepsilon)\left[
\frac{2\pi c}{3}\bigl(\sqrt2+\log(1+\sqrt2)\bigr)
-c^2+\frac{c^4}{24\pi}
\right]
+\varepsilon\rho_c(d(x,m_c)).
$$
Thus $m_c$ is the unique global minimizer.
For $x$ sufficiently close to $m_c$,
$$
F_c(x)-F_c(m_c)=\varepsilon d(x,m_c)^2,
$$
so
$$
\operatorname{Hess}^{g}F_c(m_c)(V,W)
=2\varepsilon\langle V,W\rangle,
\qquad V,W\in\T_{m_c}\M.
$$

The regular cut locus is
$$
H_{m_c}
=\{(\pi,t):-\pi<t<\pi\}
\cup\{(t,\pi):-\pi<t<\pi\},
$$
and the only exceptional cut point is $(\pi,\pi)$.
For $y=(\pi,t)$ with $-\pi<t<\pi$, we have
$$
\kappa(m_c\mid y)=\kappa(y\mid m_c)
=\frac{2\pi}{\sqrt{\pi^2+t^2}},
\qquad
n(m_c\mid y)\otimes n(m_c\mid y)=e_1\otimes e_1.
$$
For $y=(t,\pi)$ with $-\pi<t<\pi$, the corresponding
normal projection is $e_2\otimes e_2$.
This description of the cut locus is given in
\cite[Section~2.4]{HotzLeWood2024}.
Write $J_{ij}$ for the matrix entries of
$J_{\mu_\varepsilon,c}(m_c)$ in the chosen frame.
Since $c<\pi$ and the density near the cut locus is
$\psi=(1-\varepsilon)/(4\pi^2)$, we obtain
$$
J_{11}=J_{22}
=4\pi c\psi\int_{-\pi}^{\pi}
\frac{\mathrm dt}{\sqrt{\pi^2+t^2}}
=\frac{2(1-\varepsilon)c\operatorname{arsinh}(1)}{\pi},
\qquad J_{12}=J_{21}=0.
$$

Define $b(y)=a_c(\|y\|)y$ on the closed fundamental square.
This Euclidean map is Lipschitz and satisfies $-Dg_c(m_c,y)=\mathrm db_y$
for almost every $y$ in its interior, using the flat
coordinate identifications. By symmetry the integral for
the Haar component is a scalar multiple of the identity.
Moreover,
$$
\int_{(-\pi,\pi)^2}\partial_1b_1(y)\,\mathrm dy
=\int_{-\pi}^{\pi}
\bigl[b_1(\pi,t)-b_1(-\pi,t)\bigr]\,\mathrm dt
=8\pi c\operatorname{arsinh}(1).
$$
Multiplying by $(1-\varepsilon)/(4\pi^2)$ and adding the
contribution of the atom gives
$$
\Lambda_{\mu_\varepsilon,c}^{\mathrm{reg}}
=\left[
2\varepsilon
+\frac{2(1-\varepsilon)c\operatorname{arsinh}(1)}{\pi}
\right]\operatorname{Id}_2.
$$
Subtracting the cut locus correction gives
$2\varepsilon\operatorname{Id}_2$, in agreement with the
direct calculation of the population Hessian.

Write $r=\|y\|$ in the fundamental domain.
By symmetry and the fact that the score vanishes at the atom,
$$
\Sigma_c
=2(1-\varepsilon)
\mathbb E_{\mu_{\mathrm{Haar}}}\bigl[\min\{r,c\}^2\bigr]
\operatorname{Id}_2,
$$
where
$$
\mathbb E_{\mu_{\mathrm{Haar}}}\bigl[\min\{r,c\}^2\bigr]
=\frac{
\int_0^c2\pi r^3\,\mathrm dr
+c^2(4\pi^2-\pi c^2)}{4\pi^2}
=c^2-\frac{c^4}{8\pi}.
$$
The density is continuous near the cut locus.
The threshold circle $S_c(m_c)$ has zero Haar measure
and does not contain the atom, so
$\mu_\varepsilon(S_c(m_c))=0$.
Together with uniqueness and positive definiteness,
these facts yield the stated central limit theorem by
Theorem~\ref{thm:huber-clt}.

Finally, $\mathcal A_\delta(m_c)$ is the union of the
strips of width $2\delta$ around the two coordinate cut lines.
Therefore,
$$
\operatorname{vol}(\mathcal A_\delta(m_c))
=(2\pi)^2-(2\pi-2\delta)^2
=8\pi\delta-4\delta^2.
$$
The regular cut locus has total length $4\pi$ and the ratio
of the two endpoint gradient norms is $1$.
This agrees with the coefficient $8\pi$ in
Theorem~\ref{thm:sweep-volume}.
\end{proof}

\subsection{A Haar mixture on \texorpdfstring{$\mathrm{SO}(3)$}{SO(3)}}

Let $\M=\mathrm{SO}(3)$ carry the standard bi-invariant metric
normalized so that
$$
d(I,R)=\theta(R)\in[0,\pi],
$$
where $I$ is the identity rotation and $\theta(R)$ is the
rotation angle of $R$.
With this normalization, $\mathrm{SO}(3)$ is isometric to
the quotient of the three-dimensional sphere of radius $2$
by the antipodal map and has constant sectional curvature $1/4$.
Let $\mu_{\mathrm{Haar}}$ be the normalized Haar probability
measure and consider
$$
\mu_\varepsilon
=(1-\varepsilon)\mu_{\mathrm{Haar}}+\varepsilon\delta_I,
\qquad 0<\varepsilon<1.
$$
Fix $0<c<\pi$.
Choose an orthonormal basis of $\T_I\M$ and use it to
identify this tangent space isometrically with $\R^3$.
We represent linear operators by their matrices in this
basis and write $\operatorname{Id}_3$ for the identity matrix.

\begin{example}
\label{ex:so3-explicit}
The identity $I$ is the unique global minimizer of $F_c$ and
$$
\Lambda_{\mu_\varepsilon,c}
=2\varepsilon\operatorname{Id}_3,\qquad
\Sigma_c
=\frac{4(1-\varepsilon)
\bigl(3\pi c^2-2c^3-6c\cos c+6\sin c\bigr)}{9\pi}
\operatorname{Id}_3.
$$
Hence
$$
\sqrt n\,\widetilde{\operatorname{Log}}_I(\widehat m_{n,c})
\overset{d}{\longrightarrow}
N_{\T_I\M}\left(
0,
\frac{(1-\varepsilon)
\bigl(3\pi c^2-2c^3-6c\cos c+6\sin c\bigr)}
{9\pi\varepsilon^2}\operatorname{Id}_3
\right).
$$
For $0<\delta<\pi$ we also have
$$
\operatorname{vol}(\mathcal A_\delta(I))
=8\pi(\delta+\sin\delta)
=16\pi\delta-\frac{4\pi}{3}\delta^3+O(\delta^5).
$$
In particular, the first-order volume coefficient in
Theorem~\ref{thm:sweep-volume} is $16\pi$.
\end{example}

\begin{proof}
Under $\mu_{\mathrm{Haar}}$, the rotation angle has density
$(1-\cos r)/\pi$ for $0\le r\le\pi$.
Haar invariance and direct integration give
$$
F_c(R)
=(1-\varepsilon)\left[
\pi c-c^2+\frac{c^3}{3\pi}
+\frac{2(c+\sin c)}{\pi}\right]
+\varepsilon\rho_c(d(R,I)).
$$
Thus $I$ is the unique global minimizer.
For $R$ sufficiently close to $I$ we have
$F_c(R)-F_c(I)=\varepsilon d(R,I)^2$ and hence
$$
\operatorname{Hess}^{g}F_c(I)(V,W)
=2\varepsilon\langle V,W\rangle,
\qquad V,W\in\T_I\M.
$$

Under the chosen normalization,
$$
C_I=H_I=\{R:\theta(R)=\pi\}\simeq\mathbb{RP}^2.
$$
For $y\in H_I$ we have
$$
d(I,y)=\pi,
\qquad
\kappa(I\mid y)=\kappa(y\mid I)=2.
$$
The spherical quotient representation gives
$$
\operatorname{vol}(\mathrm{SO}(3))=8\pi^2,
\qquad
\mathcal H^2(H_I)=8\pi,
\qquad
\psi=\frac{1-\varepsilon}{8\pi^2},
$$
where $\psi$ is the density of $\mu_\varepsilon$ near the
cut locus.
Rotational symmetry gives
$$
\int_{H_I}n(I\mid y)\otimes n(I\mid y)
\,\mathrm d\mathcal H^2(y)
=\frac{8\pi}{3}\operatorname{Id}_3.
$$
Therefore
$$
J_{\mu_\varepsilon,c}(I)
=4c\psi\,\frac{8\pi}{3}\operatorname{Id}_3
=\frac{4(1-\varepsilon)c}{3\pi}\operatorname{Id}_3.
$$

We next compute the regular part for the Haar measure
at the same base point $I$.
The radial Hessian formula for sectional curvature $1/4$
and the conditional second moment $\operatorname{Id}_3/3$
of the unit direction given the rotation angle yield
$$
\Lambda_{\mu_{\mathrm{Haar}},c}^{\mathrm{reg}}
=\frac1{3\pi}\int_0^\pi
\left[\rho_c''(r)+\rho_c'(r)\cot\frac r2\right]
(1-\cos r)\,\mathrm dr\,\operatorname{Id}_3=\frac{4c}{3\pi}\operatorname{Id}_3.
$$
Here $\rho_c''$ is interpreted almost everywhere.
The second equality uses absolute continuity of $\rho_c'$
and the identity
$(1-\cos r)\cot(r/2)=\sin r$.
Multiplying by $1-\varepsilon$ and adding the atomic
contribution gives
$$
\Lambda_{\mu_\varepsilon,c}^{\mathrm{reg}}
=\left[2\varepsilon+\frac{4(1-\varepsilon)c}{3\pi}\right]
\operatorname{Id}_3.
$$
Subtracting the cut locus correction gives
$\Lambda_{\mu_\varepsilon,c}=2\varepsilon\operatorname{Id}_3$.
This agrees with the population Hessian computed above
through the identity
$\operatorname{Hess}^{g}F_c(I)(V,W)
=\langle\Lambda_{\mu_\varepsilon,c}V,W\rangle$.

Rotational symmetry and the vanishing score at the atom
also give
$$
\Sigma_c
=\frac{4(1-\varepsilon)}{3\pi}
\left[
\int_0^c r^2(1-\cos r)\,\mathrm dr
+c^2\int_c^\pi(1-\cos r)\,\mathrm dr
\right]\operatorname{Id}_3.
$$
Direct integration yields the stated covariance.
The density is continuous near the cut locus.
The angle density and the condition $c>0$ also imply
$\mu_\varepsilon(S_c(I))=0$.
Together with uniqueness and positive definiteness,
Theorem~\ref{thm:huber-clt} gives the stated central limit
theorem.

Finally, choose unit quaternion representatives
$p,q\in S^3$ of $I,y$ and put
$r=d(I,y)=2\arccos|\langle p,q\rangle|$.
Since $C_y$ is represented by $q^\perp\cap S^3$,
$$
d(I,C_y)
=2\arccos\sqrt{1-|\langle p,q\rangle|^2}
=\pi-r.
$$
Cut locus reciprocity therefore gives
$\mathcal A_\delta(I)=\{y:d(I,y)>\pi-\delta\}$, and hence
$$
\operatorname{vol}(\mathcal A_\delta(I))
=8\pi^2\int_{\pi-\delta}^{\pi}
\frac{1-\cos r}{\pi}\,\mathrm dr
=8\pi(\delta+\sin\delta).
$$
Its first-order coefficient is $16\pi$, which is also
twice the area of the regular cut locus and agrees with
Theorem~\ref{thm:sweep-volume}.
\end{proof}

\begin{remark}
On the circle, the correction is determined by the density
at a single cut point.
On the flat torus, it comes from the two regular cut lines.
The example on $\mathrm{SO}(3)$ exhibits a nonzero correction
in positive sectional curvature.
In the latter two mixture models, the contribution of the
uniform component to $\Lambda_{\mu_\varepsilon,c}^{\mathrm{reg}}$
is exactly canceled by the cut locus correction.
Thus the expected pointwise Hessian alone does not determine
the population Hessian.

For comparison, the cut locus of a point on the standard
sphere $S^d$ with $d\ge2$ consists only of its antipodal point.
Hence $H_m=\varnothing$ for every $m\in S^d$.
Under the assumptions of
Theorem~\ref{thm:huber-score-linearization}, we have
$$
J_{\mu,c}(m_c)=0,
\qquad
\Lambda_{\mu,c}=\Lambda_{\mu,c}^{\mathrm{reg}}.
$$
Whether a correction arises therefore cannot be determined
solely by whether the cut locus intersects the support
of the distribution.
This agrees with the phenomenon described in
\cite[Section~2.4]{HotzLeWood2024}.
\end{remark}

\section{Proof of Proposition~\ref{prop:field-difference-hotz}}
\label{sec:proof-weighted-linearization}

We prove Proposition~\ref{prop:field-difference-hotz}
by showing that the distributional covariant derivative
of $L_{\mu,w}$ has the continuous density
$A_{\mu,w}+J_{\mu,w}$.
This requires controlling the possible singularities
of $H(x\mid y)$ near the cut locus and determining how
the jumps of the logarithm map contribute to the derivative.

The required integrability, averaging, and continuity
results are established in
Lemmas~\ref{lem:weighted-adjugate},
\ref{lem:weighted-averaged-singular},
and~\ref{lem:weighted-J-continuity}
in Appendix~\ref{app:weighted-hotz}.
We combine these results below to obtain classical
$C^1$ regularity and the first-order expansion.

Fix $m,\mu,w$ as in Proposition~\ref{prop:field-difference-hotz},
and retain the notation of Section~\ref{sec:main-results}
and Appendix~\ref{app:weighted-hotz}.
Let $U$ be the normal neighborhood provided by
Lemma~\ref{lem:weighted-cut-stability}.

For linear operators $P,Q$ on a tangent space, we use the pairing
$\langle P,Q\rangle=\operatorname{tr}(P^*Q)$,
where $P^*$ denotes the adjoint with respect to the Riemannian
metric. In particular,
$$
|\langle P,Q\rangle|
\le d\|P\|_{\mathrm{op}}\|Q\|_{\mathrm{op}}.
$$
The total variation $|\mathsf S|$ of an operator-valued
Radon measure $\mathsf S$ is defined using the operator norm.

For scalar functions and vector fields, $\|\cdot\|_\infty$
denotes the supremum norm, using the absolute value and the
Riemannian norm, respectively. For a $(1,1)$-tensor field $\Phi$ on $U$,
we write
$$
\|\Phi\|_{\infty,\mathrm{op}}
=\sup_{x\in U}\|\Phi(x)\|_{\mathrm{op}}.
$$

For $d\ge2$, let $\mathsf S_y$ denote the singular operator-valued measure in
Lemma~\ref{lem:weighted-fixed-y-measure}. Let $K_0$ be the uniform bound on total variation in Lemma~\ref{lem:weighted-defect}.

\begin{proof}
We first consider the case $d\ge2$.
By Lemma~\ref{lem:weighted-cut-stability}, we have
$C_x\subset O$ and $\mu(C_x)=0$ for every $x\in U$,
which proves \textup{(i)}.

We first show that $L_{\mu,w}$ is continuous on $U$.
Fix $x\in U$ and an arbitrary sequence $x_j\to x$.
Choose a sufficiently small coordinate neighborhood of $x$
and a local smooth orthonormal frame $e_1,\ldots,e_d$ on this
neighborhood. After discarding finitely many terms, we may assume that
all $x_j$ lie in this neighborhood. It suffices to show that each component of
$L_{\mu,w}(x_j)$ in this frame converges to the corresponding
component of $L_{\mu,w}(x)$.

Since $\mu(C_x)=0$, we have $y\notin C_x$ for $\mu$-almost
every $y\in\M$. Fix such a $y$.
Cut locus reciprocity gives $x\notin C_y$.
Since $C_y$ is closed and $x_j\to x$, we also have
$x_j\notin C_y$ for all sufficiently large $j$.
Equivalently, $y\notin C_{x_j}$ for these $j$.

The map $z\mapsto\operatorname{Log}_z(y)$ is smooth on
$\M\setminus C_y$. Thus, for $k=1,\ldots,d$,
$$
\left\langle
\operatorname{Log}_{x_j}(y),e_k(x_j)
\right\rangle
\longrightarrow
\left\langle
\operatorname{Log}_{x}(y),e_k(x)
\right\rangle.
$$
So $w(y)\left\langle\operatorname{Log}_{x_j}(y),e_k(x_j)\right\rangle
\longrightarrow
w(y)\left\langle\operatorname{Log}_{x}(y),e_k(x)\right\rangle$ for $\mu$-almost every $y$.

On the other hand, since $\operatorname{diam}(\M)=D$ and $w$ is
bounded, for every $j$ and every $y\notin C_{x_j}$,
$$
\left|
w(y)
\left\langle
\operatorname{Log}_{x_j}(y),e_k(x_j)
\right\rangle
\right|
\le
|w(y)|\,\|\operatorname{Log}_{x_j}(y)\|
\le D\|w\|_\infty.
$$
For $y\in C_{x_j}$, the left-hand side is zero, so the same
bound holds for every $y\in\M$.

Since $D\|w\|_\infty$ is integrable with respect to $\mu$,
the dominated convergence theorem gives, for each $k=1,\ldots,d$,
$$
\begin{aligned}
\left\langle L_{\mu,w}(x_j),e_k(x_j)\right\rangle
&=\int_\M w(y)
\left\langle
\operatorname{Log}_{x_j}(y),e_k(x_j)
\right\rangle
\,\mathrm d\mu(y)\\
&\longrightarrow
\int_\M w(y)
\left\langle
\operatorname{Log}_{x}(y),e_k(x)
\right\rangle
\,\mathrm d\mu(y)\\
&=\left\langle L_{\mu,w}(x),e_k(x)\right\rangle.
\end{aligned}
$$
Thus $L_{\mu,w}(x_j)$ converges componentwise to $L_{\mu,w}(x)$
in the local smooth orthonormal frame.
Since the sequence $x_j\to x$ was arbitrary, $L_{\mu,w}$ is
continuous on $U$.

We next compute the distributional covariant derivative of
$L_{\mu,w}$. Let
$$
\Phi\in C_c^\infty(U;\T^*\M\otimes\T\M)
$$
be an arbitrary smooth, compactly supported test tensor.
We write $\mathsf D$ for the distributional covariant derivative
and let $\mathsf D^*$ denote the formal adjoint of $D$ with
respect to the metric-induced pairings and the Riemannian
volume measure. So
$$
\begin{aligned}
\left\langle\mathsf D L_{\mu,w},\Phi\right\rangle
&=\int_U
\left\langle L_{\mu,w}(x),\mathsf D^*\Phi(x)\right\rangle
\,\mathrm d\operatorname{vol}(x)\\
&=\int_U\int_\M w(y)
\left\langle
\operatorname{Log}_x(y),\mathsf D^*\Phi(x)
\right\rangle
\,\mathrm d\mu(y)\,\mathrm d\operatorname{vol}(x).
\end{aligned}
$$

Set $K=\operatorname{supp}\Phi\Subset U$.
For all $x,y\in\M$,
$\|w(y)\operatorname{Log}_x(y)\|\le D\|w\|_\infty$.
Moreover, $\mathsf D^*\Phi$ is smooth and supported in the compact
set $K$, so $\|\mathsf D^*\Phi\|_\infty<\infty$.
Consequently,
$$
\int_U\int_\M
\left|w(y)
\left\langle
\operatorname{Log}_x(y),\mathsf D^*\Phi(x)
\right\rangle\right|
\,\mathrm d\mu(y)\,\mathrm d\operatorname{vol}(x)\le D\|w\|_\infty
\|\mathsf D^*\Phi\|_\infty\operatorname{vol}(K)<\infty.
$$
Fubini's theorem therefore gives
$$
\begin{aligned}
\left\langle\mathsf D L_{\mu,w},\Phi\right\rangle
&=\int_\M w(y)
\left[\int_U
\left\langle
\operatorname{Log}_x(y),\mathsf D^*\Phi(x)
\right\rangle
\,\mathrm d\operatorname{vol}(x)\right]
\,\mathrm d\mu(y)\\
&=\int_\M w(y)
\left\langle
\mathsf D_x\operatorname{Log}_x(y),\Phi
\right\rangle
\,\mathrm d\mu(y).
\end{aligned}
$$

By Lemma~\ref{lem:weighted-fixed-y-measure}, for each fixed $y\in\M$,
$$
\mathsf D_x\operatorname{Log}_x(y)
=H(x\mid y)\operatorname{vol}(\mathrm dx)
+\mathsf S_y(\mathrm dx).
$$
Hence
$$
\left\langle\mathsf D L_{\mu,w},\Phi\right\rangle
=\int_\M w(y)\int_U
\left\langle H(x\mid y),\Phi(x)\right\rangle
\,\mathrm d\operatorname{vol}(x)\,\mathrm d\mu(y)+\int_\M w(y)\int_U
\langle\Phi(x),\mathrm d\mathsf S_y(x)\rangle
\,\mathrm d\mu(y).
$$

For the regular part, we estimate the tensor pairing using the
operator norm. By \eqref{eq:weighted-absolute-H} in
Lemma~\ref{lem:weighted-adjugate},
$$
\sup_{x\in U}
\int_{\M\setminus C_x}
|w(y)|\,\|H(x\mid y)\|_{\mathrm{op}}
\,\mathrm d\mu(y)<\infty.
$$
Hence
$$
\begin{aligned}
&\int_K\int_\M
|w(y)|\,
\left|\left\langle H(x\mid y),\Phi(x)\right\rangle\right|
\,\mathrm d\mu(y)\,\mathrm d\operatorname{vol}(x)
\\
&\qquad\le
d\|\Phi\|_{\infty,\mathrm{op}}
\operatorname{vol}(K)
\sup_{x\in U}
\int_{\M\setminus C_x}
|w(y)|\,\|H(x\mid y)\|_{\mathrm{op}}
\,\mathrm d\mu(y)<\infty.
\end{aligned}
$$
Fubini's theorem therefore gives
$$
\begin{aligned}
&\int_\M w(y)
\int_U\left\langle H(x\mid y),\Phi(x)\right\rangle
\,\mathrm d\operatorname{vol}(x)\,\mathrm d\mu(y)
\\
&\qquad=
\int_U\left\langle
\int_{\M\setminus C_x}w(y)H(x\mid y)\,\mathrm d\mu(y),
\Phi(x)\right\rangle\,\mathrm d\operatorname{vol}(x)
\\
&\qquad=
\int_U\left\langle A_{\mu,w}(x),\Phi(x)\right\rangle
\,\mathrm d\operatorname{vol}(x).
\end{aligned}
$$

For the singular part, Lemma~\ref{lem:weighted-defect} yields
$|\mathsf S_y|(\M)\le K_0$ for every $y\in\M$. Thus
$$
\int_\M |w(y)|\,|\mathsf S_y|(U)\,\mathrm d\mu(y)
\le\|w\|_\infty K_0<\infty.
$$
By the kernel measurability argument in
Lemma~\ref{lem:weighted-averaged-singular}, each component of
$y\mapsto w(y)\mathsf S_y|_U$ in a local smooth orthonormal frame
is a Borel signed measure kernel. The preceding total variation
bound ensures that the integrals of these kernels are well defined
and allows us to apply the Tonelli--Fubini theorem componentwise.
By Lemma~\ref{lem:weighted-averaged-singular}, the integral
defining $J_{\mu,w}(x)$ converges absolutely for every $x\in U$,
and we have
$$
\int_\M w(y)\mathsf S_y|_U\,\mathrm d\mu(y)
=J_{\mu,w}\,\operatorname{vol}|_U.
$$
Consequently,
$$
\int_\M w(y)
\int_U\langle\Phi(x),\mathrm d\mathsf S_y(x)\rangle
\,\mathrm d\mu(y)
=\int_U\left\langle J_{\mu,w}(x),\Phi(x)\right\rangle
\,\mathrm d\operatorname{vol}(x).
$$

Combining the regular and singular parts, we obtain
$$
\left\langle\mathsf D L_{\mu,w},\Phi\right\rangle
=\int_U\left\langle
A_{\mu,w}(x)+J_{\mu,w}(x),\Phi(x)
\right\rangle\,\mathrm d\operatorname{vol}(x).
$$
Since this holds for every
$\Phi\in C_c^\infty(U;\T^*\M\otimes\T\M)$,
we obtain
$$
\mathsf D L_{\mu,w}
=\bigl(A_{\mu,w}+J_{\mu,w}\bigr)\operatorname{vol}|_U.
$$

The fields $A_{\mu,w}$ and $J_{\mu,w}$ are continuous on $U$
by Lemmas~\ref{lem:weighted-adjugate}
and~\ref{lem:weighted-J-continuity}.
We now use the distributional identity above to establish
classical $C^1$ regularity.

Fix a coordinate neighborhood $V\Subset U$ with local
coordinates $(x^1,\ldots,x^d)$.
Write $\partial_i=\partial/\partial x^i$ and let
$\Gamma^j_{ik}$ denote the Christoffel symbols in these
coordinates. Repeated upper and lower indices are summed
from $1$ to $d$.
We write
$$
L_{\mu,w}=L^j\partial_j,
\qquad
(A_{\mu,w}+J_{\mu,w})(\partial_i)
=(A_{\mu,w}+J_{\mu,w})^j{}_i\partial_j.
$$
The distributional identity becomes
$$
\partial_iL^j
=(A_{\mu,w}+J_{\mu,w})^j{}_i-\Gamma^j_{ik}L^k
$$
in these coordinates.
The right-hand side is continuous because
$A_{\mu,w}$, $J_{\mu,w}$, and $L_{\mu,w}$ are continuous
and the Christoffel symbols are smooth.

By \cite[Theorem~6.10]{LiebLoss2001}, each $L^j$ belongs
to $C^1(V)$ and the same identity holds for its classical
partial derivatives. Hence
$$
\bigl(D_{\partial_i}L_{\mu,w}\bigr)^j
=\partial_iL^j+\Gamma^j_{ik}L^k
=(A_{\mu,w}+J_{\mu,w})^j{}_i.
$$
Since $V\Subset U$ was arbitrary, we obtain
$L_{\mu,w}\in C^1(U;\T\M)$ and
$$
DL_{\mu,w}(x)=A_{\mu,w}(x)+J_{\mu,w}(x),
\qquad x\in U.
$$
This proves \textup{(ii)}.

For $x$ sufficiently close to $m$, integrating the covariant
derivative along the unique minimizing geodesic from $m$
to $x$ gives
\begin{equation}
\label{eq:weighted-frechet-bound}
\begin{aligned}
&\left\|
\Pi_{x,m}L_{\mu,w}(x)-L_{\mu,w}(m)
-\bigl[A_{\mu,w}(m)+J_{\mu,w}(m)\bigr]
\operatorname{Log}_m(x)
\right\|\\
&\qquad\le d(x,m)
\sup_{d(z,m)\le d(x,m)}
\left\|
\Pi_{z,m}DL_{\mu,w}(z)\Pi_{m,z}
-DL_{\mu,w}(m)
\right\|_{\mathrm{op}}\\
&\qquad=o(d(x,m)).
\end{aligned}
\end{equation}
The last equality follows from continuity of $DL_{\mu,w}$.
This proves \textup{(iii)}.

We now consider the case $d=1$.
The manifold $\M$ is isometric to a circle of circumference
$L$. In arc-length coordinates, we represent tangent vectors
and linear operators by their scalar components.
Let $a(x)$ be the antipodal point of $x$, so that
$H_x=C_x=\{a(x)\}$.

For fixed $y$, the scalar representation of
$x\mapsto\operatorname{Log}_x(y)$ has slope $-1$ and an
upward jump of size $L$ at $a(y)$. Hence
$$
\mathsf D_x\operatorname{Log}_x(y)
=-\operatorname{vol}(\mathrm dx)
+L\delta_{a(y)}(\mathrm dx).
$$
For every $x\in U$, the measure $\mu$ has density $\psi$
near $a(x)$. Since $a$ is an arc-length-preserving involution,
averaging the preceding identity with weight $w$ gives
$$
A_{\mu,w}(x)=-\int_\M w\,\mathrm d\mu,
\qquad
J_{\mu,w}(x)=Lw(a(x))\psi(a(x)).
$$
The latter agrees with \eqref{eq:weighted-J-definition}
because $\mathcal H^0$ is the counting measure and
$$
d(x,a(x))=\frac L2,
\qquad
\kappa(x\mid a(x))=\kappa(a(x)\mid x)=2.
$$

Both fields are continuous and their defining integrals
converge absolutely.
The same continuity and weak derivative arguments as above
give $L_{\mu,w}\in C^1(U;\T\M)$ with
$DL_{\mu,w}=A_{\mu,w}+J_{\mu,w}$.
The one-dimensional fundamental theorem of calculus then
gives \eqref{eq:weighted-frechet-bound}.
\end{proof}

\section{Proofs of Theorem~\ref{thm:huber-score-linearization}
and Corollary~\ref{cor:huber-risk-quadratic}}
\label{sec:proof-linearization}

This section uses the notation and zero-extension conventions
of Section~\ref{subsec:notation}. For $y\notin C_x\cup\{x\}$, write
$$
\nu_x(y)=\frac{\operatorname{Log}_x(y)}{d(x,y)}.
$$
This is the unit initial velocity of the unique minimizing
geodesic from $x$ to $y$. For $r>0$ with $r\ne c$, we have
$$
a_c'(r)=
\begin{cases}
0,&0<r<c,\\[1mm]
-\dfrac{2c}{r^2},&r>c,
\end{cases}
\qquad |ra_c'(r)|\le2.
$$

\subsection{Two auxiliary lemmas}
\label{subsec:huber-decomposition}

Throughout this subsection, set $w(y)=a_c(d(m_c,y))$ and
let $L_{\mu,w}$ be the weighted logarithm vector field in
Proposition~\ref{prop:field-difference-hotz}.
For $x$ sufficiently close to $m_c$, define
$$
R_c^{\mathrm{rad}}(x)
=\int_\M
\bigl[a_c(d(x,y))-a_c(d(m_c,y))\bigr]
\Pi_{x,m_c}\operatorname{Log}_x(y)\,\mathrm d\mu(y).
$$
Since $L_{\mu,w}(m_c)=G_c(m_c)$, we have
$$
\Pi_{x,m_c}G_c(x)-G_c(m_c)
=\Pi_{x,m_c}L_{\mu,w}(x)-L_{\mu,w}(m_c)
+R_c^{\mathrm{rad}}(x).
$$
The following lemmas give first-order expansions of the two
terms on the right.

\begin{lemma}
\label{lem:weighted-hotz-linearization}
Under assumptions \textup{(i)} and \textup{(iii)} of
Theorem~\ref{thm:huber-score-linearization}, the integral
defining $J_{\mu,c}(m_c)$ converges absolutely.
As $x\to m_c$,
$$
\begin{aligned}
\Pi_{x,m_c}L_{\mu,w}(x)-L_{\mu,w}(m_c)
&=
\left[\int_{\M\setminus C_{m_c}}
a_c(d(m_c,y))H(m_c\mid y)\,\mathrm d\mu(y)\right]
\operatorname{Log}_{m_c}(x)\\
&\qquad+J_{\mu,c}(m_c)\operatorname{Log}_{m_c}(x)
+o(d(x,m_c)).
\end{aligned}
$$
\end{lemma}

\begin{proof}
The weight $w$ is continuous on $\M$.
Apply Proposition~\ref{prop:field-difference-hotz} with
$m=m_c$ and use the definition of $A_{\mu,w}$ and
$J_{\mu,w}(m_c)=J_{\mu,c}(m_c)$.
\end{proof}

\begin{lemma}
\label{lem:huber-radial-weight-variation}
Under assumptions \textup{(i)}, \textup{(iii)}, and
\textup{(iv)} of Theorem~\ref{thm:huber-score-linearization},
define the linear operator
$$
B_{\mu,c}(m_c)
=2c\int_{\{y\in\M\setminus C_{m_c}:\,d(m_c,y)>c\}}
\frac{\nu_{m_c}(y)\otimes\nu_{m_c}(y)}{d(m_c,y)}
\,\mathrm d\mu(y).
$$
Then
$$
R_c^{\mathrm{rad}}(x)
=B_{\mu,c}(m_c)\operatorname{Log}_{m_c}(x)
+o(d(x,m_c))
\qquad\text{as }x\to m_c.
$$
If $c\ge D$, the integration domain is empty and
$B_{\mu,c}(m_c)=0$.
\end{lemma}

\begin{proof}
If $c\ge D$, then $a_c\equiv2$ on $[0,D]$, so
$R_c^{\mathrm{rad}}=0$ and $B_{\mu,c}(m_c)=0$.
The conclusion follows immediately.

Suppose that $0<c<D$. Assumptions \textup{(iii)} and
\textup{(iv)} of Theorem~\ref{thm:huber-score-linearization}
and Lemma~\ref{lem:weighted-cut-stability} give
$$
\mu(C_{m_c}\cup S_c(m_c))=0.
$$
The integrand defining $B_{\mu,c}(m_c)$ has operator norm
at most $2$ because $2c/d(m_c,y)\le2$ on its domain and
$\|\nu_{m_c}(y)\otimes\nu_{m_c}(y)\|_{\mathrm{op}}=1$.
Thus this operator is well defined.

In the expressions below, the operator
$$
\frac{2c}{r}\mathbf1_{\{r>c\}}
\bigl(\nu_{m_c}(y)\otimes\nu_{m_c}(y)\bigr),
\qquad r=d(m_c,y),
$$
is set to zero on $C_{m_c}\cup S_c(m_c)\cup\{m_c\}$.
Define
$$
\begin{aligned}
\varepsilon_c(x,y)
={}&\bigl[a_c(d(x,y))-a_c(d(m_c,y))\bigr]
\Pi_{x,m_c}\operatorname{Log}_x(y)\\
&-\frac{2c}{r}\mathbf1_{\{r>c\}}
\bigl(\nu_{m_c}(y)\otimes\nu_{m_c}(y)\bigr)
\operatorname{Log}_{m_c}(x).
\end{aligned}
$$
For fixed $y\notin C_{m_c}\cup S_c(m_c)\cup\{m_c\}$,
cut locus reciprocity and closedness of $C_y$ imply that
$x\mapsto d(x,y)$ and $x\mapsto\operatorname{Log}_x(y)$
are smooth near $m_c$. The first variation formula and smoothness of the logarithm
map give
$$
\begin{aligned}
a_c(d(x,y))-a_c(r)
&=-a_c'(r)
\langle\nu_{m_c}(y),\operatorname{Log}_{m_c}(x)\rangle
+o(d(x,m_c)),\\
\Pi_{x,m_c}\operatorname{Log}_x(y)
&=r\nu_{m_c}(y)+o(1).
\end{aligned}
$$
Multiplying these expansions and using
$-ra_c'(r)=(2c/r)\mathbf1_{\{r>c\}}$ yields
$$
\begin{aligned}
&[a_c(d(x,y))-a_c(r)]
\Pi_{x,m_c}\operatorname{Log}_x(y)\\
&\qquad=
\frac{2c}{r}\mathbf1_{\{r>c\}}
\bigl(\nu_{m_c}(y)\otimes\nu_{m_c}(y)\bigr)
\operatorname{Log}_{m_c}(x)
+o(d(x,m_c)).
\end{aligned}
$$
By the definition of $\varepsilon_c$, we obtain
$$
\frac{\|\varepsilon_c(x,y)\|}{d(x,m_c)}
\longrightarrow0
\qquad\text{as }x\to m_c.
$$
When $y=m_c$, the weight difference vanishes for
$d(x,m_c)<c$, so $\varepsilon_c(x,m_c)=0$ there.
The preceding convergence thus holds for $\mu$-almost
every $y$.

Since $\operatorname{Lip}(a_c)\le2/c$ and
$\|\Pi_{x,m_c}\operatorname{Log}_x(y)\|\le D$,
the operator norm bound above yields
$$
\frac{\|\varepsilon_c(x,y)\|}{d(x,m_c)}
\le\frac{2D}{c}+2.
$$
Dominated convergence now gives
$$
\frac{
\|R_c^{\mathrm{rad}}(x)
-B_{\mu,c}(m_c)\operatorname{Log}_{m_c}(x)\|
}{d(x,m_c)}\le
\int_\M\frac{\|\varepsilon_c(x,y)\|}{d(x,m_c)}
\,\mathrm d\mu(y)
\longrightarrow0.
$$
This proves Lemma~\ref{lem:huber-radial-weight-variation}.
\end{proof}

\subsection{Proof of Theorem~\ref{thm:huber-score-linearization}}

\begin{proof}
Lemma~\ref{lem:weighted-cut-stability} gives
$\mu(C_{m_c})=0$.
If $0<c<D$, assumption \textup{(iv)} gives
$\mu(S_c(m_c))=0$. If $c>D$, this set is empty.
For $c=D$, we have $S_D(m_c)\subset C_{m_c}$:
otherwise the distance function $z\mapsto d(m_c,z)$ would
attain its maximum at a smooth point $y\ne m_c$ where its
gradient has norm $1$, a contradiction.
Thus $\mu(C_{m_c}\cup S_c(m_c))=0$ in every case.

We first establish differentiability of $F_c$ at $m_c$.
For small $h\in\T_{m_c}\M$, set
$$
q_{m_c,y}(h)
=\rho_c(d(\operatorname{Exp}_{m_c}h,y))
-\rho_c(d(m_c,y))
+\langle g_c(m_c,y),h\rangle.
$$
The first variation formula gives
$q_{m_c,y}(h)=o(\|h\|)$ for
$y\notin C_{m_c}\cup\{m_c\}$.
For $y=m_c$, we have $q_{m_c,m_c}(h)=\|h\|^2$ when $h$
is sufficiently small.
The $2c$-Lipschitz property of $\rho_c$ and the bound
$\|g_c(m_c,y)\|\le2c$ give
$$
\frac{|q_{m_c,y}(h)|}{\|h\|}\le4c.
$$
Since $\mu(C_{m_c})=0$, dominated convergence yields
$$
\int_\M\frac{|q_{m_c,y}(h)|}{\|h\|}
\,\mathrm d\mu(y)\longrightarrow0.
$$
Hence
$$
F_c(\operatorname{Exp}_{m_c}h)-F_c(m_c)
=-\langle G_c(m_c),h\rangle+o(\|h\|).
$$
Thus $F_c$ is Fr\'echet differentiable at $m_c$ with
$\operatorname{grad}F_c(m_c)=-G_c(m_c)$.
Assumption \textup{(ii)} now gives $G_c(m_c)=0$.

For $y\notin C_{m_c}\cup S_c(m_c)\cup\{m_c\}$, put
$r=d(m_c,y)$. Covariant differentiation gives
$$
D_Vg_c(m_c,y)
=a_c(r)H(m_c\mid y)V
-ra_c'(r)
\bigl(\nu_{m_c}(y)\otimes\nu_{m_c}(y)\bigr)V.
$$
At $y=m_c$, the derivative is $D_Vg_c(m_c,m_c)=-2V$.
Proposition~\ref{prop:field-difference-hotz} with $w\equiv1$
shows that $H(m_c\mid\cdot)$ is absolutely integrable.
Together with $0<a_c\le2$ and $|ra_c'(r)|\le2$, this proves
absolute convergence of the integral defining
$\Lambda_{\mu,c}^{\mathrm{reg}}$.
Since $-ra_c'(r)=\frac{2c}{r}\mathbf1_{\{r>c\}}$
for $r>0$ with $r\ne c$, we obtain
$$
-\Lambda_{\mu,c}^{\mathrm{reg}}
=\int_{\M\setminus C_{m_c}}
a_c(d(m_c,y))H(m_c\mid y)\,\mathrm d\mu(y)
+B_{\mu,c}(m_c).
$$
Here we used $\mu(C_{m_c}\cup S_c(m_c))=0$ and the
definition of $B_{\mu,c}(m_c)$.

Lemma~\ref{lem:weighted-hotz-linearization} ensures absolute
convergence of the integral defining $J_{\mu,c}(m_c)$.
Combining the two auxiliary lemmas with the decomposition
above gives
$$
\begin{aligned}
\Pi_{x,m_c}G_c(x)-G_c(m_c)
&=\left[
\int_{\M\setminus C_{m_c}}
a_c(d(m_c,y))H(m_c\mid y)\,\mathrm d\mu(y)
\right]\operatorname{Log}_{m_c}(x)\\
&\qquad+
\bigl[B_{\mu,c}(m_c)+J_{\mu,c}(m_c)\bigr]
\operatorname{Log}_{m_c}(x)+o(d(x,m_c))\\
&=\bigl[-\Lambda_{\mu,c}^{\mathrm{reg}}+J_{\mu,c}(m_c)\bigr]
\operatorname{Log}_{m_c}(x)+o(d(x,m_c))\\
&=-\Lambda_{\mu,c}\operatorname{Log}_{m_c}(x)
+o(d(x,m_c)).
\end{aligned}
$$
Since $G_c(m_c)=0$, this is the asserted expansion.
If $H_{m_c}=\varnothing$, then $J_{\mu,c}(m_c)=0$ by its
definition.

It remains to verify self-adjointness.
For $y\notin C_{m_c}$ and $V,W\in\T_{m_c}\M$,
$$
\langle H(m_c\mid y)V,W\rangle
=-\operatorname{Hess}^{g}
\left(\frac12d(\,\cdot\,,y)^2\right)(m_c)(V,W).
$$
The symmetry of the Hessian shows that $H(m_c\mid y)$
is self-adjoint.
The operators $\nu_{m_c}(y)\otimes\nu_{m_c}(y)$ and
$n(m_c\mid y)\otimes n(m_c\mid y)$ are also self-adjoint.
The derivative formula and absolute convergence of the
integrals therefore imply that
$\Lambda_{\mu,c}^{\mathrm{reg}}$, $J_{\mu,c}(m_c)$, and
$\Lambda_{\mu,c}$ are self-adjoint.
\end{proof}

\subsection{Proof of Corollary~\ref{cor:huber-risk-quadratic}}

\begin{proof}
By Lemma~\ref{lem:weighted-cut-stability}, choose a normal
neighborhood $U$ of $m_c$ such that $\mu(C_z)=0$ for every
$z\in U$. The differentiability argument in the preceding
proof applies at each $z\in U$ and gives
$$
F_c(\operatorname{Exp}_z h)-F_c(z)
=-\langle G_c(z),h\rangle+o(\|h\|).
$$
Thus $F_c$ is differentiable on $U$ with
$\operatorname{grad}F_c=-G_c$.

To prove continuity of $G_c$, let $z_j\to z$ in $U$.
Choose a local smooth orthonormal frame near $z$ and discard
finitely many terms so that it is defined at every $z_j$.
All convergence statements in this paragraph refer to
components in this frame.
For $y\notin C_z\cup\{z\}$, cut locus reciprocity and
closedness of $C_y$ imply that
$x\mapsto\operatorname{Log}_x(y)$ is smooth near $z$.
Continuity of $a_c$ therefore gives
$g_c(z_j,y)\to g_c(z,y)$.
For $y=z$, we have
$$
g_c(z_j,z)=2\operatorname{Log}_{z_j}(z)
\longrightarrow0=g_c(z,z)
$$
for all sufficiently large $j$.
Since $\mu(C_z)=0$ and $\|g_c(z_j,y)\|\le2c$,
dominated convergence gives $G_c(z_j)\to G_c(z)$.
Hence $F_c\in C^1(U)$.

Theorem~\ref{thm:huber-score-linearization} shows that
$\operatorname{grad}F_c$ is Fr\'echet differentiable at $m_c$
with
$$
D_V\operatorname{grad}F_c(m_c)=\Lambda_{\mu,c}V,
\qquad V\in\T_{m_c}\M.
$$
Consequently,
$$
\operatorname{Hess}^{g}F_c(m_c)(V,W)
=\langle\Lambda_{\mu,c}V,W\rangle,
\qquad V,W\in\T_{m_c}\M.
$$

For a sufficiently small nonzero $u\in\T_{m_c}\M$, write $v=u/\|u\|$ and set
$$
\gamma(s)=\operatorname{Exp}_{m_c}(sv),
\qquad 0\le s\le \|u\|.
$$
The curve lies in $U$ when $u$ is small.
By the score expansion,
$$
\Pi_{\gamma(s),m_c}G_c(\gamma(s))
=-s\Lambda_{\mu,c}v+o(s),
$$
where the remainder is $o(s)$ uniformly over unit vectors
$v$. Since
$\Pi_{\gamma(s),m_c}\dot\gamma(s)=v$ and parallel transport
preserves inner products,
$$
\frac{\mathrm d}{\mathrm ds}F_c(\gamma(s))
=-\langle G_c(\gamma(s)),\dot\gamma(s)\rangle
=s\langle v,\Lambda_{\mu,c}v\rangle+o(s).
$$
Integrating gives
$$
\begin{aligned}
F_c(\operatorname{Exp}_{m_c}u)-F_c(m_c)
&=\int_0^{\|u\|}
\bigl[s\langle v,\Lambda_{\mu,c}v\rangle+o(s)\bigr]
\,\mathrm ds\\
&=\frac12\langle u,\Lambda_{\mu,c}u\rangle
+o(\|u\|^2).
\end{aligned}
$$

If $\Lambda_{\mu,c}$ is strictly positive definite, let
$\lambda_{\min}>0$ be its smallest eigenvalue.
The expansion implies that for all sufficiently small $u$,
$$
F_c(\operatorname{Exp}_{m_c}u)-F_c(m_c)
\ge\frac{\lambda_{\min}}2\|u\|^2
-\frac{\lambda_{\min}}4\|u\|^2
=\frac{\lambda_{\min}}4\|u\|^2.
$$
Taking $a=\lambda_{\min}/4$ and using
$\|\operatorname{Log}_{m_c}(x)\|=d(x,m_c)$ gives the
asserted lower bound for $d(x,m_c)<\delta$ with sufficiently
small $\delta>0$.
\end{proof}

\section{Proof of Theorem~\ref{thm:huber-clt}}
\label{sec:proof-clt}

For every $\mu$-integrable function $f$, write
$$
Pf=\int_\M f(y)\,\mathrm d\mu(y),
\qquad
P_nf=\frac1n\sum_{i=1}^n f(X_i).
$$

The $2c$-Lipschitz property of the Huber loss and the triangle
inequality give
$$
|\rho_c(d(x,y))-\rho_c(d(z,y))|
\le2c\,d(x,z),
\qquad x,z,y\in\M.
$$
\subsection{Strong consistency}
\begin{lemma}
\label{lem:huber-selection-consistency}
Fix $c\in(0,\infty)$. Suppose that $\M$ is compact and
$m_c$ is the unique global minimizer of $F_c$.
Then every measurable selection $\widehat m_{n,c}\in E_{n,c}$
satisfies
$$
\widehat m_{n,c}\overset{\mathrm{a.s.}}{\longrightarrow}m_c.
$$
\end{lemma}

\begin{proof}
For each integer $k\ge1$, compactness gives a finite
$1/k$-net $x_{k,1},\ldots,x_{k,N_k}$ in $\M$.
Fix $x\in\M$ and choose a point $x_{k,j}$ that satisfies $d(x,x_{k,j})<1/k$.
The triangle inequality and the preceding Lipschitz bound give
$$
\begin{aligned}
|F_{n,c}(x)-F_c(x)|
&\le |F_{n,c}(x)-F_{n,c}(x_{k,j})|
+|F_{n,c}(x_{k,j})-F_c(x_{k,j})|\\
&\quad+|F_c(x_{k,j})-F_c(x)|\\
&\le4c\,d(x,x_{k,j})
+|F_{n,c}(x_{k,j})-F_c(x_{k,j})|\\
&\le\frac{4c}{k}
+\max_{1\le j\le N_k}
|F_{n,c}(x_{k,j})-F_c(x_{k,j})|.
\end{aligned}
$$
Since $0\le\rho_c(d(x_{k,j},X_i))\le\rho_c(D)$, the strong
law of large numbers gives
$$
\max_{1\le j\le N_k}|F_{n,c}(x_{k,j})-F_c(x_{k,j})|
\longrightarrow0
\qquad\text{almost surely for every }k\ge1.
$$
These convergences hold simultaneously on an event of
probability one. Taking the limit superior and then letting
$k\to\infty$ yields
$$
\sup_{x\in\M}|F_{n,c}(x)-F_c(x)|
\overset{\mathrm{a.s.}}{\longrightarrow}0.
$$

Work on this event. For every sample minimizer,
$$
\begin{aligned}
0\le F_c(\widehat m_{n,c})-F_c(m_c)
&\le |F_c(\widehat m_{n,c})-F_{n,c}(\widehat m_{n,c})|
   +|F_{n,c}(m_c)-F_c(m_c)|\\
&\le2\sup_{x\in\M}|F_{n,c}(x)-F_c(x)|
\longrightarrow0.
\end{aligned}
$$
For $\varepsilon>0$, the set
$A_\varepsilon=\{x\in\M:d(x,m_c)\ge\varepsilon\}$ is
compact. If it is nonempty, continuity of $F_c$ and uniqueness
of its minimizer give
$$
\inf_{x\in A_\varepsilon}\{F_c(x)-F_c(m_c)\}>0.
$$
Thus every sample minimizer eventually lies outside
$A_\varepsilon$. If $A_\varepsilon$ is empty, this is immediate.
Since $\varepsilon>0$ is arbitrary, the claimed convergence
holds on the same event of probability one.
\end{proof}

\subsection{Proof of Theorem~\ref{thm:huber-clt}}

\begin{proof}
Fix an orthonormal basis of $\T_{m_c}\M$ and use the
corresponding coordinates to identify this space with $\R^d$.
Retain $U_0$ and $r_0$ from Section~\ref{subsec:huber-clt}.
Set $\Theta=B_{r_0}(0)$ and define
$$
q_u(y)=-\rho_c(d(\operatorname{Exp}_{m_c}u,y)),
\qquad u\in\Theta,
\qquad
\widehat u_n
=\widetilde{\operatorname{Log}}_{m_c}(\widehat m_{n,c}).
$$
Then
$$
Pq_u=-F_c(\operatorname{Exp}_{m_c}u),
\qquad
P_nq_u=-F_{n,c}(\operatorname{Exp}_{m_c}u).
$$
The function $(u,y)\mapsto q_u(y)$ is continuous and
$|q_u(y)|\le\rho_c(D)$. The random vector $\widehat u_n$
is measurable and takes values in $\Theta$.
By Lemma~\ref{lem:huber-selection-consistency},
$$
\widehat u_n\longrightarrow0\quad\text{almost surely},
\qquad
\mathbb P(\widehat m_{n,c}\in U_0)\longrightarrow1.
$$
We verify the remaining assumptions of
\cite[Theorem~5.23]{vdV1998}.

\textup{(i)} By Theorem~\ref{thm:huber-score-linearization},
$\mu(C_{m_c})=0$.
For $y\notin C_{m_c}\cup\{m_c\}$, the distance to $y$ is
smooth near $m_c$ and $\rho_c$ is continuously differentiable.
Hence $u\mapsto q_u(y)$ is differentiable at $0$ with
$$
\left.\frac{\mathrm d}{\mathrm dt}q_{th}(y)\right|_{t=0}
=\frac{\rho_c'(d(m_c,y))}{d(m_c,y)}
\langle\operatorname{Log}_{m_c}(y),h\rangle
=\langle g_c(m_c,y),h\rangle,
\qquad h\in\T_{m_c}\M.
$$
For $y=m_c$, we have $q_u(m_c)=-\|u\|^2$ for sufficiently
small $u$, so the derivative at $0$ is zero.
Thus $u\mapsto q_u(y)$ is Fr\'echet differentiable at $0$
for $\mu$-almost every $y$ with derivative vector
$g_c(m_c,y)$ in the chosen coordinates.

\textup{(ii)} Smoothness of the exponential map and compactness
of $\overline\Theta$ give
$$
L=\sup_{z\in\overline\Theta}
\|(\mathrm d\operatorname{Exp}_{m_c})_z\|_{\mathrm{op}}
<\infty.
$$
For $u,v\in\Theta$, the curve
$t\mapsto\operatorname{Exp}_{m_c}((1-t)u+tv)$ with
$0\le t\le1$ lies in $U_0$ because $\Theta$ is convex.
Its length gives
$$
d(\operatorname{Exp}_{m_c}u,\operatorname{Exp}_{m_c}v)
\le\int_0^1
\|(\mathrm d\operatorname{Exp}_{m_c})_{(1-t)u+tv}(v-u)\|
\,\mathrm dt
\le L\|u-v\|.
$$
Consequently,
$$
|q_u(y)-q_v(y)|
\le2c\,d(\operatorname{Exp}_{m_c}u,\operatorname{Exp}_{m_c}v)
\le2cL\|u-v\|,
\qquad y\in\M.
$$
The constant function $2cL$ has finite second moment $4c^2L^2<\infty$.

\textup{(iii)} Uniqueness of $m_c$ and injectivity of
$\operatorname{Exp}_{m_c}$ on $\Theta$ imply that $0$ is
the unique maximizer of $u\mapsto Pq_u$ on $\Theta$.
By Corollary~\ref{cor:huber-risk-quadratic},
$$
Pq_u-Pq_0
=-\frac12\langle u,\Lambda_{\mu,c}u\rangle
+o(\|u\|^2)
\qquad\text{as }u\to0.
$$
The matrix representing $-\Lambda_{\mu,c}$ is symmetric
and negative definite, hence nonsingular.

\textup{(iv)} Define the maximization error by
$$
\Delta_n=\sup_{u\in\Theta}P_nq_u-P_nq_{\widehat u_n}\ge0.
$$
For every sample, $u\mapsto P_nq_u$ is continuous, so
$$
\sup_{u\in\Theta}P_nq_u
=\sup_{u\in\Theta\cap\mathbb Q^d}P_nq_u.
$$
The latter is a supremum of countably many measurable random
variables and is therefore measurable. Joint continuity of
$q_u(y)$ and measurability of $\widehat u_n$ also imply that
$P_nq_{\widehat u_n}$ is measurable. Thus $\Delta_n$ is a
finite measurable random variable.

On the event $\{\widehat m_{n,c}\in U_0\}$, we have
$\operatorname{Exp}_{m_c}(\widehat u_n)=\widehat m_{n,c}$.
Since $\widehat m_{n,c}$ is a global sample minimizer,
$$
P_nq_{\widehat u_n}
=-F_{n,c}(\widehat m_{n,c})
\ge-F_{n,c}(\operatorname{Exp}_{m_c}u)
=P_nq_u,
\qquad u\in\Theta.
$$
Hence $\Delta_n=0$ on this event. For every $\varepsilon>0$,
$$
\mathbb P(n\Delta_n>\varepsilon)
\le\mathbb P(\widehat m_{n,c}\notin U_0)
\longrightarrow0.
$$
Thus $\Delta_n=o_p(n^{-1})$, which gives
$$
P_nq_{\widehat u_n}
=\sup_{u\in\Theta}P_nq_u-o_p(n^{-1}).
$$

The theorem now yields
$$
\sqrt n\,\widehat u_n
=\Lambda_{\mu,c}^{-1}
\frac1{\sqrt n}\sum_{i=1}^n g_c(m_c,X_i)+o_p(1).
$$
Since $\mathbb E[g_c(m_c,X_1)]=G_c(m_c)=0$ and
$\|g_c(m_c,X_1)\|\le2c$, the multivariate central limit
theorem and Slutsky's theorem give
$$
\sqrt n\,\widetilde{\operatorname{Log}}_{m_c}(\widehat m_{n,c})
\overset{d}{\longrightarrow}
N_{\T_{m_c}\M}\left(
0,\Lambda_{\mu,c}^{-1}\Sigma_c\Lambda_{\mu,c}^{-1}
\right).
$$
\end{proof}

\newpage
\appendix
\section{Geometric and measure-theoretic tools for differentiating weighted logarithm fields}
\label{app:weighted-hotz}
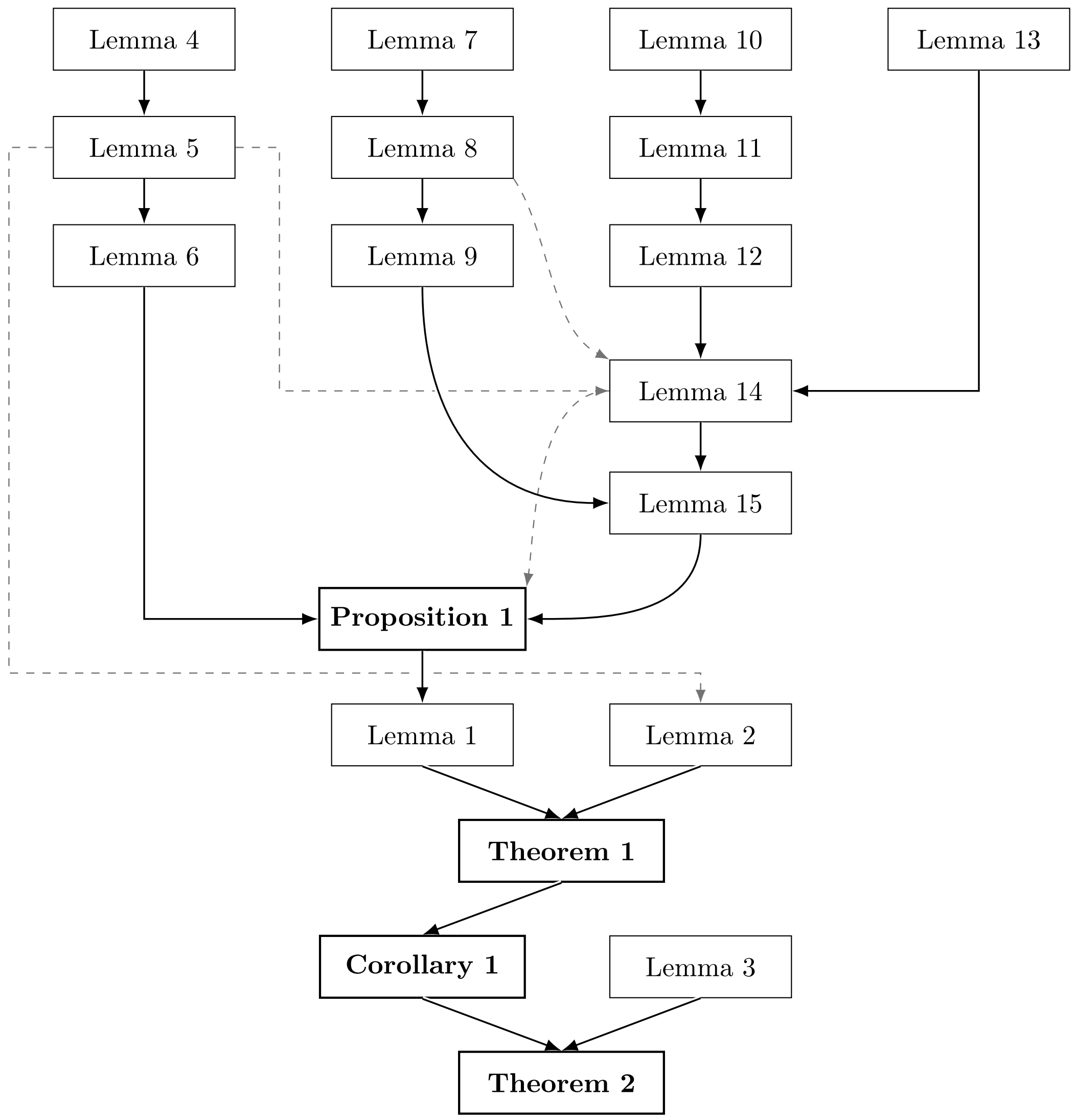
\begin{figure}[H]
\centering
\begin{tikzpicture}[
  >=Latex,
  box/.style={draw=black,rectangle,align=center,
    inner xsep=4pt,inner ysep=5pt,font=\small,
    minimum height=8mm,minimum width=2.35cm},
  main/.style={box,line width=.8pt,font=\small\bfseries,
    minimum width=2.65cm},
  arr/.style={->,line width=.65pt,
    preaction={draw=white,line width=2pt,solid}},
  auxarr/.style={->,dashed,line width=.45pt,draw=black!55}
]
\path[use as bounding box]  (-7.3,-14) rectangle (7.3,0.5);
\node[box] (l4) at (-5.4,0) {Lemma~\ref{lem:cut-time-closedness}};
\node[box] (l5) at (-5.4,-1.4) {Lemma~\ref{lem:weighted-cut-stability}};
\node[box] (l6) at (-5.4,-2.8) {Lemma~\ref{lem:weighted-adjugate}};
\node[box] (l7) at (-1.8,0) {Lemma~\ref{lem:weighted-fixed-y-measure}};
\node[box] (l8) at (-1.8,-1.4) {Lemma~\ref{lem:weighted-defect}};
\node[box] (l9) at (-1.8,-2.8) {Lemma~\ref{lem:weighted-defect-continuity}};
\node[box] (l10) at (1.8,0) {Lemma~\ref{lem:weighted-incidence-structure}};
\node[box] (l11) at (1.8,-1.4) {Lemma~\ref{lem:weighted-incidence-jacobian}};
\node[box] (l12) at (1.8,-2.8) {Lemma~\ref{lem:weighted-incidence-coarea}};
\node[box] (l13) at (5.4,0) {Lemma~\ref{lem:weighted-angle-comparison}};
\node[box] (l14) at (1.8,-4.55) {Lemma~\ref{lem:weighted-averaged-singular}};
\node[box] (l15) at (1.8,-6) {Lemma~\ref{lem:weighted-J-continuity}};
\node[main] (p1) at (-1.8,-7.5) {Proposition~\ref{prop:field-difference-hotz}};
\node[box] (l1) at (-1.8,-9) {Lemma~\ref{lem:weighted-hotz-linearization}};
\node[box] (l2) at (1.8,-9) {Lemma~\ref{lem:huber-radial-weight-variation}};
\node[main] (t1) at (0,-10.5) {Theorem~\ref{thm:huber-score-linearization}};
\node[main] (c1) at (-1.8,-12) {Corollary~\ref{cor:huber-risk-quadratic}};
\node[box] (l3) at (1.8,-12) {Lemma~\ref{lem:huber-selection-consistency}};
\node[main] (t2) at (0,-13.5) {Theorem~\ref{thm:huber-clt}};
\draw[auxarr] (l5.east) -- (-3.65,-1.4) -- (-3.65,-4.55) -- (l14.west);
\draw[auxarr] (l8.south east) .. controls (-0.15,-2.5) and (-0.15,-3.75) .. (l14.north west);
\draw[auxarr] (l14.west) .. controls (-0.35,-4.55) and (-0.35,-6.4) .. (p1.north east);
\draw[auxarr] (l5.west) -- (-7.15,-1.4) -- (-7.15,-8.2) -- (1.8,-8.2) -- (l2.north);
\draw[arr] (l4.south) -- (l5.north);
\draw[arr] (l5.south) -- (l6.north);
\draw[arr] (l7.south) -- (l8.north);
\draw[arr] (l8.south) -- (l9.north);
\draw[arr] (l10.south) -- (l11.north);
\draw[arr] (l11.south) -- (l12.north);
\draw[arr] (l12.south) -- (l14.north);
\draw[arr] (l14.south) -- (l15.north);
\draw[arr] (l13.south) -- (5.4,-4.55) -- (l14.east);
\draw[arr] (l9.south) .. controls (-1.8,-4.6) and (-1.25,-6) .. (l15.west);
\draw[arr] (l6.south) -- (-5.4,-7.5) -- (p1.west);
\draw[arr] (l15.south) .. controls (1.8,-7.5) and (0.5,-7.5) .. (p1.east);
\draw[arr] (p1.south) -- (l1.north);
\draw[arr] (l1.south) -- (t1.north);
\draw[arr] (l2.south) -- (t1.north);
\draw[arr] (t1.south) -- (c1.north);
\draw[arr] (c1.south) -- (t2.north);
\draw[arr] (l3.south) -- (t2.north);
\end{tikzpicture}
\caption{Main proof dependencies. Arrows point from a result
to a result that uses it. Solid arrows show the main chains
and dashed arrows show supporting uses.
Some transitive dependencies are omitted.}
\label{fig:proof-dependency}
\end{figure}

Figure~\ref{fig:proof-dependency} summarizes the main
proof dependencies. The proof of
Theorem~\ref{thm:sweep-volume} is detailed separately in
Figure~\ref{fig:sweep-proof-dependency}.
We retain the assumptions and notation of
Proposition~\ref{prop:field-difference-hotz}, together with
the conventions of Sections~\ref{subsec:notation}
and~\ref{sec:proof-weighted-linearization}.
For $p\in\M$, let $d_p(y)=d(p,y)$.
\subsection{Localization and the regular part}
\label{subsec:weighted-localization}
\begin{lemma}
\label{lem:cut-time-closedness}
Let
$$
S\M=\{(x,v):x\in\M,\ v\in\T_x\M,\ \|v\|=1\}
$$
denote the unit tangent bundle.
For $(x,v)\in S\M$, define the cut time by
$$
b(x,v)=\sup\left\{
s>0:
d\bigl(x,\operatorname{Exp}_x(tv)\bigr)=t, \forall t\in[0,s]
\right\}.
$$
Then the following statements hold.

\textup{(i)} We have $0<b(x,v)\le D$, and the map
$b:S\M\to(0,D]$ is continuous. Moreover, the set
$$
\mathcal C=\{(x,y)\in\M\times\M:y\in C_x\}
$$
is compact.

\textup{(ii)} For every $m\in\M$ and $\delta>0$, the set
$\mathcal A_\delta(m)$ is Borel.

\textup{(iii)} For every $m\in\M$ and every open set $V$ containing
$C_m$, there exists $\delta_V>0$ such that
$$
\mathcal A_\delta(m)\subset V,
\qquad 0<\delta<\delta_V.
$$
\end{lemma}

\begin{proof}
By \cite[Theorem~10.33]{Lee2018}, the cut-time function
$b:S\M\to(0,\infty]$ is continuous.
Since $\M$ is compact, $0<b(x,v)\le D$.
The continuous map
$$
S\M\longrightarrow\M\times\M,
\qquad
(x,v)\longmapsto
\left(x,\operatorname{Exp}_x\bigl(b(x,v)v\bigr)\right)
$$
has image $\mathcal C$.
Since $S\M$ is compact, so is $\mathcal C$.
This proves \textup{(i)}.

Fix $m\in\M$ and $\delta>0$, and choose a sequence
$0<r_j\uparrow\delta$ with $r_j<\delta$ for every $j$.
Then
$$
\mathcal A_\delta(m)
=\bigcup_j\bigcup_{x\in\overline B_{r_j}(m)}C_x.
$$
For each $j$, the set $\displaystyle\bigcup_{x\in\overline B_{r_j}(m)}C_x$
is the projection onto the second factor of the compact set
$\mathcal C\cap\bigl(\overline B_{r_j}(m)\times\M\bigr)$, 
and is therefore compact.
Thus $\mathcal A_\delta(m)$ is a countable union of compact sets,
and hence is Borel. This proves \textup{(ii)}.

Fix $m\in\M$ and an open set $V$ containing $C_m$.
If the conclusion were false, then for every $j\ge1$ we could choose
$x_j\in B_{1/j}(m)$ and $y_j\in C_{x_j}\setminus V$.
Thus $x_j\to m$.

By compactness of $\M$, there is a subsequence such that
$y_{k_j}\to y\in\M$.
Since $x_{k_j}\to m$, we have
$$
(x_{k_j},y_{k_j})\longrightarrow(m,y).
$$
Each $(x_{k_j},y_{k_j})$ belongs to $\mathcal C$, and $\mathcal C$
is closed. Hence $(m,y)\in\mathcal C$, or equivalently $y\in C_m$.
On the other hand, $y_{k_j}\in\M\setminus V$, and $\M\setminus V$
is closed, so $y\in\M\setminus V$.
Consequently, $y\in C_m\setminus V$, contradicting $C_m\subset V$.
This proves \textup{(iii)}.
\end{proof}
\begin{lemma}
\label{lem:weighted-cut-stability}
Choose open sets $W_0,W_1$ such that
$$
C_m\subset W_0\Subset W_1\Subset O.
$$
Then there exists a normal neighborhood $U$ of $m$ such that
\begin{equation}
\label{eq:weighted-localization}
C_x\subset W_0,\qquad \mu(C_x)=0,
\qquad x\in\overline U.
\end{equation}
\end{lemma}

\begin{proof}
Open sets $W_0,W_1$ exist because $C_m$ is compact.
By Lemma~\ref{lem:cut-time-closedness}\textup{(iii)},
there exists an open neighborhood $V$ of $m$ such that
$C_x\subset W_0$ for every $x\in V$.
Choose a sufficiently small normal ball $U$ centered at $m$
with $\overline U\subset V$.

For every $x\in\overline U$, we have $C_x\subset W_0\subset O$.
Since $C_x$ has zero Riemannian volume
by \cite[Theorem~10.34(a)]{Lee2018}
and $\mu|_O$ is absolutely continuous with respect to
$\operatorname{vol}|_O$, it follows that $\mu(C_x)=0$.
\end{proof}
For the subsequent arguments concerning weighted integrals,
fix the normal neighborhood $U$ and the open sets $W_0,W_1$
from Lemma~\ref{lem:weighted-cut-stability}, satisfying
\eqref{eq:weighted-localization}.
Choose $\chi\in C^\infty(\M)$ such that
$$
0\le\chi\le1,\qquad
\chi|_{W_0}=1,\qquad
\operatorname{supp}\chi\subset W_1,
$$
and set
$$
M_\psi=\sup_{\overline{W_1}}\psi<\infty.
$$
These choices of $U,W_0,W_1,\chi$ and $M_\psi$ remain fixed
throughout the subsequent arguments concerning weighted
integrals in this appendix.
\begin{lemma}
\label{lem:weighted-adjugate}
Under the assumptions of Proposition~\ref{prop:field-difference-hotz},
\begin{equation}
\label{eq:weighted-absolute-H}
\sup_{x\in U}
\int_{\M\setminus C_x}
|w(y)|\,\|H(x\mid y)\|_{\mathrm{op}}
\,\mathrm d\mu(y)
<\infty.
\end{equation}
Moreover, $A_{\mu,w}$ is continuous on $U$.
In particular, taking $w\equiv1$ shows that
$H(m\mid\cdot)$ is absolutely integrable with respect to $\mu$.
\end{lemma}
\begin{proof}

For $x\in\M$ and $v\in\T_x\M$, define linear maps
$$
V_{x,v},B_{x,v}:\T_x\M
\longrightarrow\T_{\operatorname{Exp}_x(v)}\M
$$
by
$$
V_{x,v}a
=\left.\frac{\mathrm d}{\mathrm dt}
\operatorname{Exp}_x(v+ta)\right|_{t=0},\qquad B_{x,v}a=\left.\frac{\mathrm d}{\mathrm dt}
\operatorname{Exp}_{\gamma(t)}(v(t))\right|_{t=0},
\qquad a\in\T_x\M,
$$
where $\gamma$ is any smooth curve with $\gamma(0)=x$ and
$\dot\gamma(0)=a$, and $v(t)$ is parallel along $\gamma$
with $v(0)=v$. Let
$$
P^{x,v}_{\operatorname{Exp}_x(v),x}:
\T_{\operatorname{Exp}_x(v)}\M\longrightarrow\T_x\M
$$
denote parallel transport from $\operatorname{Exp}_x(v)$ to $x$
along the geodesic $s\mapsto\operatorname{Exp}_x(sv)$,
$0\le s\le1$, and set
$$
\mathcal V_{x,v}
=P^{x,v}_{\operatorname{Exp}_x(v),x}V_{x,v},\qquad
\mathcal B_{x,v}
=P^{x,v}_{\operatorname{Exp}_x(v),x}B_{x,v}.
$$
Both operators depend smoothly on $(x,v)$.

For $y\notin C_x$, set $v=\operatorname{Log}_x(y)$.
Differentiating $\operatorname{Exp}_z(\operatorname{Log}_z y)=y$
at $z=x$ in the direction $a\in\T_x\M$ gives
$$
B_{x,v}a+V_{x,v}H(x\mid y)a=0.
$$
Since $V_{x,v}$ is invertible on this domain, parallel
transport back to $\T_x\M$ yields
$$
H(x\mid y)=-\mathcal V_{x,v}^{-1}\mathcal B_{x,v}.
$$
Since $\mathcal V_{x,0}=\operatorname{Id}$, continuity along
radial segments gives $\det\mathcal V_{x,v}>0$.
Thus the volume Jacobian of $\operatorname{Exp}_x$ is
$\det(\mathcal V_{x,v})$, and
\begin{equation}
\label{eq:weighted-adjugate-cancellation}
\det(\mathcal V_{x,v})\,H(x\mid\operatorname{Exp}_x(v))
=-\operatorname{adj}(\mathcal V_{x,v})\mathcal B_{x,v}.
\end{equation}
Here $\operatorname{adj}$ denotes the classical adjugate.
Since the adjugate depends polynomially on the matrix entries,
the right-hand side of
\eqref{eq:weighted-adjugate-cancellation} is continuous on
the compact disk bundle $\{(x,v)\in\T\M:\|v\|\le D\}$
and is therefore uniformly bounded.

Using $\chi$, split the integral in
\eqref{eq:weighted-A-definition} into
$$
A_{\mu,w}(x)
=A_{\mu,w}^{\mathrm{near}}(x)+A_{\mu,w}^{\mathrm{far}}(x),
$$
where
$$
\begin{aligned}
A_{\mu,w}^{\mathrm{near}}(x)
&=\int_{\M\setminus C_x}
\chi(y)w(y)H(x\mid y)\,\mathrm d\mu(y),\\
A_{\mu,w}^{\mathrm{far}}(x)
&=\int_{\M\setminus C_x}
(1-\chi(y))w(y)H(x\mid y)\,\mathrm d\mu(y).
\end{aligned}
$$
We establish absolute convergence and continuity of both terms.

Since $C_x\subset W_0$ for every $x\in\overline U$ and
$\chi\equiv1$ on $W_0$,
$$
\bigl(\overline U\times\operatorname{supp}(1-\chi)\bigr)
\cap\mathcal C=\varnothing.
$$
Hence $H$ is smooth on a neighborhood of the compact set
$\overline U\times\operatorname{supp}(1-\chi)$
and is uniformly bounded in operator norm on this set. Since $w$ is bounded and
$\mu$ is finite, the far term satisfies
$$
\sup_{x\in U}\int_{\M\setminus C_x}
(1-\chi(y))|w(y)|\,\|H(x\mid y)\|_{\mathrm{op}}
\,\mathrm d\mu(y)<\infty,
$$
and $A_{\mu,w}^{\mathrm{far}}$ is continuous by dominated convergence.

For the near term, define
$$
\varphi(y)=
\begin{cases}
w(y)\chi(y)\psi(y),&y\in O,\\
0,&y\notin O.
\end{cases}
$$
Since $\operatorname{supp}\chi\subset W_1\Subset O$,
$\varphi$ is continuous and bounded on $\M$.
For $x\in U$, define
$$
I_x=\left\{
t\theta:\theta\in\T_m\M,\ \|\theta\|=1,
    \ 0\le t<b(x,\Pi_{m,x}\theta)
\right\}.
$$
The change of variables $y=\operatorname{Exp}_x(\Pi_{m,x}z), z\in I_x$, together with $\mu|_O=\psi\,\operatorname{vol}|_O$ and
\eqref{eq:weighted-adjugate-cancellation}, yields a constant $C<\infty$ (independent of $x\in U$) such that
$$
\int_{\M\setminus C_x}
\chi(y)|w(y)|\,\|H(x\mid y)\|_{\mathrm{op}}
\,\mathrm d\mu(y)\le C\int_{\|z\|\le D}\mathbf1_{I_x}(z)\,\mathrm dz<\infty.
$$
Here $\mathrm dz$ denotes the Lebesgue measure on $\T_m\M$
induced by $g_m$.
Combining this estimate with the far-term bound proves
\eqref{eq:weighted-absolute-H}.

The same change of variables gives
\begin{equation}
\label{eq:weighted-A-fixed-domain}
\Pi_{x,m}A_{\mu,w}^{\mathrm{near}}(x)\Pi_{m,x}
=\int_{\|z\|\le D}\mathbf1_{I_x}(z)T(x,z)\,\mathrm dz,
\end{equation}
where
$$
T(x,z)
=-\varphi\bigl(\operatorname{Exp}_x(\Pi_{m,x}z)\bigr)
\Pi_{x,m}
\operatorname{adj}(\mathcal V_{x,\Pi_{m,x}z})
\mathcal B_{x,\Pi_{m,x}z}\Pi_{m,x}.
$$
The map $T$ is continuous on
$U\times\{z\in\T_m\M:\|z\|\le D\}$ and uniformly bounded
there, since parallel transport is an isometry and the
adjugate expression above is uniformly bounded.

Let $x_j\to x$ in $U$.
For each $\theta\in\T_m\M$ with $\|\theta\|=1$,
Lemma~\ref{lem:cut-time-closedness}\textup{(i)} gives
$$
b(x_j,\Pi_{m,x_j}\theta)
\longrightarrow b(x,\Pi_{m,x}\theta).
$$
For $z=t\theta$ with $t\ge0$, if
$t\ne b(x,\Pi_{m,x}\theta)$, the defining inequality
for membership in $I_{x_j}$ agrees with that for $I_x$
for all sufficiently large $j$. Hence $\mathbf1_{I_{x_j}}(z)\longrightarrow\mathbf1_{I_x}(z)$ except possibly on the set
$$
\left\{
b(x,\Pi_{m,x}\theta)\theta:
\theta\in\T_m\M,\ \|\theta\|=1
\right\}.
$$
This set meets each radial ray in at most one point and
therefore has Lebesgue measure zero by polar coordinates
and Fubini's theorem.
Dominated convergence applied to
\eqref{eq:weighted-A-fixed-domain} now proves the continuity
of $A_{\mu,w}^{\mathrm{near}}$.
Together with the far term, this proves the continuity of
$A_{\mu,w}$ on $U$.
\end{proof}
\subsection{Distributional derivatives and measure estimates}
\label{subsec:weighted-distributional-measures}

For the remainder of this appendix, assume $d\ge2$.

\begin{lemma}
\label{lem:weighted-fixed-y-measure}
Fix $y\in\M$. The distributional covariant derivative
$\mathsf D_x\operatorname{Log}_x(y)$ is an operator-valued
Radon measure satisfying
\begin{equation}
\label{eq:weighted-fixed-y-identity}
\mathsf D_x\operatorname{Log}_x(y)
=H(x\mid y)\,\operatorname{vol}(\mathrm dx)
+\mathsf S_y(\mathrm dx),
\end{equation}
where
\begin{equation}
\label{eq:weighted-Sy}
\mathsf S_y(\mathrm dx)
=d(x,y)\kappa(x\mid y)
\bigl(n(x\mid y)\otimes n(x\mid y)\bigr)
\,\mathcal H^{d-1}|_{H_y}(\mathrm dx)
\end{equation}
is positive semidefinite. In particular,
$|\mathsf S_y|(Q_y)=0$, and
$\mathsf D_x\operatorname{Log}_x(y)$ has no Cantor part.
\end{lemma}

\begin{proof}
We write $\operatorname{Hess}$ for the distributional Hessian.
By \cite[Proposition~2.6]{MMU2014},
$$
\operatorname{grad}d_y
\in\mathrm{SBV}_{\mathrm{loc}}(\M\setminus\{y\}).
$$
Here $\mathrm{SBV}_{\mathrm{loc}}$ denotes the class of vector
fields of locally bounded variation whose distributional
derivatives have no Cantor part.
Moreover, \cite[Theorem~2.10]{MMU2014} gives, for any local
smooth vector fields $X,Y$,
$$
\operatorname{Hess}d_y(X,Y)
=\operatorname{Hess}^g d_y(X,Y)\,\operatorname{vol}-\kappa(\cdot\mid y)
\langle n(\cdot\mid y),X\rangle
\langle n(\cdot\mid y),Y\rangle
\,\mathcal H^{d-1}|_{H_y}.
$$

Work in a coordinate neighborhood compactly contained in
$\M\setminus\{y\}$.
Since $d_y$ is locally Lipschitz, the BV chain and product
rules give
$$
d_y\operatorname{grad}d_y
\in\mathrm{SBV}_{\mathrm{loc}}(\M\setminus\{y\})
$$
and, for local smooth vector fields $X,Y$,
$$
\begin{aligned}
\operatorname{Hess}\left(\frac12d_y^2\right)(X,Y)
={}&\mathrm d d_y(X)\,\mathrm d d_y(Y)\,\operatorname{vol}
+d_y\operatorname{Hess}d_y(X,Y)\\
={}&-\langle H(\cdot\mid y)X,Y\rangle\,\operatorname{vol}\\
&-d_y\kappa(\cdot\mid y)
\langle n(\cdot\mid y),X\rangle
\langle n(\cdot\mid y),Y\rangle
\,\mathcal H^{d-1}|_{H_y}.
\end{aligned}
$$
For the absolutely continuous term, we used the identity
$$
\operatorname{Log}_x(y)
=-\operatorname{grad}\left(\frac12d_y^2\right)(x),
\qquad x\notin C_y,
$$
and the definition of $H$.
Since $C_y$ has zero Riemannian volume, this identity also
holds almost everywhere on $\M$.
Metric compatibility therefore yields the scalar distributional identity
$$
\left\langle
\bigl(\mathsf D_x\operatorname{Log}_x(y)\bigr)X,Y
\right\rangle
=-\operatorname{Hess}\left(\frac12d_y^2\right)(X,Y).
$$
Comparing these identities for arbitrary $X,Y$ proves
\eqref{eq:weighted-fixed-y-identity} on $\M\setminus\{y\}$,
with $\mathsf S_y$ given by \eqref{eq:weighted-Sy}.

Since $\frac12d_y^2$ is smooth near $y$,
\eqref{eq:weighted-fixed-y-identity} holds on all of $\M$.
The $\mathrm{SBV}_{\mathrm{loc}}$ property above and smoothness
near $y$ also show that
$\mathsf D_x\operatorname{Log}_x(y)$ has no Cantor part.

Finally, for any local smooth vector field $X$,
$$
\langle\mathsf S_yX,X\rangle
=d_y\kappa(\cdot\mid y)
\langle n(\cdot\mid y),X\rangle^2
\,\mathcal H^{d-1}|_{H_y}\ge0.
$$
Thus $\mathsf S_y$ is positive semidefinite, and
$|\mathsf S_y|(Q_y)=0$ follows directly from
\eqref{eq:weighted-Sy}.
\end{proof}
We use the sign convention
$\Delta=\operatorname{div}\operatorname{grad}$
for the Laplace--Beltrami operator and interpret
$\Delta(\frac12d_x^2)$ in the distributional sense.

\begin{lemma}
\label{lem:weighted-defect}
There exists a constant $C_0\ge1$ depending only on $(\M,g)$
such that
$$
\lambda_x
=C_0d\,\operatorname{vol}
-\Delta\left(\frac12d_x^2\right)
$$
is a positive Radon measure for every $x\in\M$.
Its total mass is
$\lambda_x(\M)=C_0d\operatorname{vol}(\M)=:K_0$.
For every Borel set $E\subset\M$ we have
$$
\int_{H_x\cap E}
d(x,y)\kappa(y\mid x)\,\mathrm d\mathcal H^{d-1}(y)
\le\lambda_x(E).
$$
Moreover, $\lambda_x(Q_x)=0, |\mathsf S_x|(\M)\le K_0$.
\end{lemma}

\begin{proof}
Choose $a>0$ such that the sectional curvature of $\M$ is
bounded below by $-a^2$, and set
$$
C_0=\max_{0\le r\le D}\max\{1,ar\coth(ar)\},
$$
where $ar\coth(ar)$ is assigned its limiting value $1$ at
$r=0$.
By the Hessian comparison theorem
\cite[Theorem~1.1]{MMU2014}, for any local smooth vector
field $X$ on $\M\setminus(C_x\cup\{x\})$,
$$
\operatorname{Hess}^{g}d_x(X,X)
\le a\coth(ad_x)
\left(\|X\|^2-\bigl(\mathrm d d_x(X)\bigr)^2\right).
$$
Since $\|\operatorname{grad}d_x\|=1$ there, it follows that
$$
\begin{aligned}
\operatorname{Hess}^{g}\left(\frac12d_x^2\right)(X,X)
&=\bigl(\mathrm d d_x(X)\bigr)^2
+d_x\operatorname{Hess}^{g}d_x(X,X)\\
&\le\bigl(\mathrm d d_x(X)\bigr)^2
+ad_x\coth(ad_x)
\left(\|X\|^2-\bigl(\mathrm d d_x(X)\bigr)^2\right)\\
&\le C_0\|X\|^2.
\end{aligned}
$$
The final bound also holds at $x$ by smoothness of
$\frac12d_x^2$ near $x$.
Using the definition of $H$ and taking the trace gives
$$
C_0d+\operatorname{tr}H(y\mid x)\ge0,
\qquad y\notin C_x.
$$

Apply Lemma~\ref{lem:weighted-fixed-y-measure} with $x$ and
$y$ interchanged. For a local smooth orthonormal frame
$e_1,\ldots,e_d$, metric compatibility gives
$$
\operatorname{Hess}\left(\frac12d_x^2\right)(e_i,e_i)
=-\langle H(\cdot\mid x)e_i,e_i\rangle\,\operatorname{vol}
-\langle\mathsf S_xe_i,e_i\rangle.
$$
Summing over $i$ yields the scalar measure identity
$$
\Delta\left(\frac12d_x^2\right)
=-\operatorname{tr}H(\cdot\mid x)\,\operatorname{vol}
-d_x\kappa(\cdot\mid x)\,\mathcal H^{d-1}|_{H_x}.
$$
Consequently,
$$
\lambda_x
=\bigl(C_0d+\operatorname{tr}H(\cdot\mid x)\bigr)
\operatorname{vol}
+d_x\kappa(\cdot\mid x)\,\mathcal H^{d-1}|_{H_x}.
$$
Both terms are nonnegative, so $\lambda_x$ is a positive
Radon measure.

Since $\M$ is compact without boundary, testing the
distributional Laplacian against $1$ gives
$$
\left\langle\Delta\left(\frac12d_x^2\right),1\right\rangle
=\int_\M\frac12d_x^2\,\Delta1\,\mathrm d\operatorname{vol}
=0.
$$
Hence $\lambda_x(\M)=C_0d\operatorname{vol}(\M)=K_0$.
For every Borel set $E\subset\M$, we have
$$
\begin{aligned}
\lambda_x(E)
&=\int_E
\bigl(C_0d+\operatorname{tr}H(y\mid x)\bigr)
\,\mathrm d\operatorname{vol}(y)+\int_{H_x\cap E}
d(x,y)\kappa(y\mid x)\,\mathrm d\mathcal H^{d-1}(y)\\
&\ge\int_{H_x\cap E}
d(x,y)\kappa(y\mid x)\,\mathrm d\mathcal H^{d-1}(y).
\end{aligned}
$$
Both integrals vanish for $E=Q_x$ because
$\operatorname{vol}(Q_x)=0$ and $Q_x\cap H_x=\varnothing$.
Thus $\lambda_x(Q_x)=0$.

Finally, $\|n(y\mid x)\otimes n(y\mid x)\|_{\mathrm{op}}=1$
on $H_x$. Thus \eqref{eq:weighted-Sy}, with $x$ and $y$
interchanged, gives
$$
|\mathsf S_x|(\M)
=\int_{H_x}d(x,y)\kappa(y\mid x)
\,\mathrm d\mathcal H^{d-1}(y)
\le\lambda_x(\M)=K_0.
$$
\end{proof}

\begin{lemma}
\label{lem:weighted-defect-continuity}
If $x_j\to x$ in $\M$, then $\lambda_{x_j}$ converges
weakly to $\lambda_x$.
\end{lemma}

\begin{proof}
The estimate
$$
\sup_{y\in\M}
\left|\frac12d(x_j,y)^2-\frac12d(x,y)^2\right|
\le D\,d(x_j,x)
$$
shows that $\frac12d_{x_j}^2$ converges uniformly to
$\frac12d_x^2$ on $\M$.
For any $\zeta\in C^\infty(\M)$, the definition of the
distributional Laplacian gives
$$
\int_\M\zeta\,\mathrm d\lambda_x
=C_0d\int_\M\zeta\,\mathrm d\operatorname{vol}
-\int_\M\frac12d_x^2\,\Delta\zeta
\,\mathrm d\operatorname{vol}.
$$
Uniform convergence therefore implies
$$
\int_\M\zeta\,\mathrm d\lambda_{x_j}
\longrightarrow\int_\M\zeta\,\mathrm d\lambda_x.
$$
By Lemma~\ref{lem:weighted-defect},
$$
\lambda_{x_j}(\M)=\lambda_x(\M)=K_0.
$$
Since $C^\infty(\M)$ is uniformly dense in $C(\M)$,
this common mass bound extends the convergence to every
$\zeta\in C(\M)$, proving the claim.
\end{proof}
\subsection{Incidence geometry and coarea}
\label{subsec:weighted-incidence-geometry}
\begin{lemma}
\label{lem:weighted-incidence-structure}
Let
$$
\Sigma_H=\{(x,y)\in\M\times\M:y\in H_x\}.
$$
Every point of $\Sigma_H$ has an open neighborhood
$\mathcal O\subset\M\times\M$ on which the distance is the
minimum of two smooth geodesic length branches
$\varphi_1$ and $\varphi_2$. Their difference
$h=\varphi_1-\varphi_2$ satisfies
$$
\Sigma_H\cap\mathcal O
=\mathcal C\cap\mathcal O
=\{(x,y)\in\mathcal O:h(x,y)=0\}.
$$
On this zero set we have
$$
\|\operatorname{grad}_x h\|=\kappa(x\mid y),
\qquad
\|\operatorname{grad}_y h\|=\kappa(y\mid x).
$$

Consequently, $\Sigma_H$ is a smooth hypersurface of
$\M\times\M$ and is relatively open in $\mathcal C$.
The coordinate projections $\pi_1(x,y)=x$ and $\pi_2(x,y)=y$
restrict to smooth submersions from $\Sigma_H$ to $\M$.
\end{lemma}

\begin{proof}
Fix $(x_*,y_*)\in\Sigma_H$. Parametrize the two minimizing
geodesics from $x_*$ to $y_*$ at constant speed on $[0,1]$
and denote their initial velocities by $v_1^*,v_2^*\in\T_{x_*}\M$.
Since the endpoints are nonconjugate along both geodesics,
the differential of $\operatorname{Exp}_{x_*}$ is invertible
at each $v_i^*$.

We claim that there is an open neighborhood $\mathcal O$ of
$(x_*,y_*)$ in $\M\times\M$ on which the two geodesics extend
to smooth branches with nonzero initial velocities and
nonconjugate endpoints. Moreover, every minimizing geodesic
from $x$ to $y$ for $(x,y)\in\mathcal O$ belongs to one of
these branches.

Apply the implicit function theorem to
$\operatorname{Exp}_x(v)=y$ with $(x,y)$ as parameters.
This gives disjoint open neighborhoods $N_i$ of $(x_*,v_i^*)$
in $\T\M$ such that for each nearby pair $(x,y)$ there is a
unique $v_i(x,y)\in\T_x\M$ satisfying
$$
(x,v_i(x,y))\in N_i,
\qquad
\operatorname{Exp}_x(v_i(x,y))=y.
$$
The two solutions depend smoothly on $(x,y)$. Since $v_i^*\ne0$
and the differential of $\operatorname{Exp}_{x_*}$ is
invertible at $v_i^*$, the initial velocities remain nonzero
and the endpoints remain nonconjugate for all nearby pairs.

To prove the assertion about minimizing geodesics, suppose
there were pairs $(x_j,y_j)\to(x_*,y_*)$ with minimizing
initial velocities $w_j$ satisfying $(x_j,w_j)\notin N_1\cup N_2$.
Since $\|w_j\|=d(x_j,y_j)\le D$, compactness of the closed
disk bundle of radius $D$ gives a subsequence along which
$(x_j,w_j)\to(x_*,w_*)$. By continuity,
$$
\operatorname{Exp}_{x_*}(w_*)=y_*,
\qquad
\|w_*\|=d(x_*,y_*).
$$
Thus $w_*$ is a minimizing initial velocity from $x_*$ to $y_*$
and must equal $v_1^*$ or $v_2^*$. Since $N_1\cup N_2$ is an open neighborhood of the limit, the subsequence eventually lies in it. This contradicts our
choice of the sequence, for which $(x_j,w_j)\notin N_1\cup N_2$ for every $j$.
This proves the claim.

For $(x,y)\in\mathcal O$ the geodesic
$t\mapsto\operatorname{Exp}_x(tv_i(x,y))$ on $[0,1]$
has length
$$
\varphi_i(x,y)=\|v_i(x,y)\|,
\qquad i=1,2.
$$
It follows that $d(x,y)=\min\{\varphi_1(x,y),\varphi_2(x,y)\}$ on $\mathcal O$.
If the lengths differ, the shorter branch gives the unique
minimizing geodesic. Its endpoints are nonconjugate, so
$y\notin C_x$. If the lengths agree, the two branches give
exactly two minimizing geodesics with nonconjugate endpoints,
so $y\in H_x$. Hence $h=\varphi_1-\varphi_2$ satisfies
$$
\Sigma_H\cap\mathcal O
=\mathcal C\cap\mathcal O
=\{(x,y)\in\mathcal O:h(x,y)=0\}.
$$

At a point $(x,y)$ in this zero set, let $u_i\in\T_x\M$ and
$z_i\in\T_y\M$ denote the unit initial and terminal tangent
vectors of the $i$th geodesic. The first variation formula
for length gives
$$
(\mathrm d\varphi_i)_{(x,y)}(\xi,\eta)
=-\langle u_i,\xi\rangle+\langle z_i,\eta\rangle,
\qquad (\xi,\eta)\in\T_x\M\times\T_y\M.
$$
Subtracting the two identities yields
$$
\operatorname{grad}_x h=u_2-u_1,
\qquad
\operatorname{grad}_y h=z_1-z_2.
$$
For fixed $x$, the functions $\varphi_i(x,\cdot)$ are the
two local distance branches at $y$. Reversing the geodesics
gives the corresponding branches at $x$ for fixed $y$.
The definition of $\kappa$ now gives the asserted norms of
the two partial gradients. Both are nonzero on the zero set,
so $\mathrm dh\ne0$ there. The implicit function theorem
shows that $\Sigma_H\cap\mathcal O$ is a smooth hypersurface.
Its equality with $\mathcal C\cap\mathcal O$ also shows that
$\Sigma_H$ is relatively open in $\mathcal C$.

Finally, at $(x,y)\in\Sigma_H\cap\mathcal O$ the tangent space is
$$
\T_{(x,y)}\Sigma_H
=\bigl\{(\xi,\eta)\in\T_x\M\times\T_y\M:
\langle\operatorname{grad}_x h,\xi\rangle
+\langle\operatorname{grad}_y h,\eta\rangle=0\bigr\}.
$$
Given $\xi\in\T_x\M$, take
$$
\eta
=-\frac{\langle\operatorname{grad}_x h,\xi\rangle}
        {\|\operatorname{grad}_y h\|^2}
  \operatorname{grad}_y h.
$$
Then $(\xi,\eta)\in\T_{(x,y)}\Sigma_H$ and
$\mathrm d\pi_1(\xi,\eta)=\xi$. Thus the differential of
$\pi_1|_{\Sigma_H}$ is surjective. Interchanging the two
variables proves the same assertion for $\pi_2|_{\Sigma_H}$.
\end{proof}
Equip $\Sigma_H$ with the metric induced by the product metric
on $\M\times\M$. For $i=1,2$, let $\mathcal J_d\pi_i$ denote
the $d$-dimensional Jacobian of $\pi_i:\Sigma_H\to\M$.

\begin{lemma}
\label{lem:weighted-incidence-jacobian}
For every $(x,y)\in\Sigma_H$ we have
$$
\mathcal J_d\pi_1(x,y)=\frac{\kappa(y\mid x)}{\sqrt{\kappa(x\mid y)^2+\kappa(y\mid x)^2}},\qquad\mathcal J_d\pi_2(x,y)=\frac{\kappa(x\mid y)}{\sqrt{\kappa(x\mid y)^2+\kappa(y\mid x)^2}}.
$$
In particular,
$$
\frac{\mathcal J_d\pi_2(x,y)}{\mathcal J_d\pi_1(x,y)}
=\frac{\kappa(x\mid y)}{\kappa(y\mid x)}.
$$
\end{lemma}

\begin{proof}
Fix $(x,y)\in\Sigma_H$ and choose a local defining function
$h$ as in Lemma~\ref{lem:weighted-incidence-structure}.
At $(x,y)$ we have
$$
\|\operatorname{grad}_x h\|=\kappa(x\mid y)>0,
\qquad
\|\operatorname{grad}_y h\|=\kappa(y\mid x)>0.
$$
The tangent space $\T_{(x,y)}\Sigma_H=\ker\mathrm dh_{(x,y)}$
contains the mutually orthogonal subspaces
$$
(\operatorname{grad}_x h)^\perp\times\{0\},
\qquad
\{0\}\times(\operatorname{grad}_y h)^\perp.
$$
Choose a unit tangent vector $(\xi,\eta)$ so that $\T_{(x,y)}\Sigma_H$ is the orthogonal direct sum of these two subspaces and $\operatorname{span}\{(\xi,\eta)\}$. Then $\xi$ is parallel to $\operatorname{grad}_x h$ and $\eta$ is parallel to $\operatorname{grad}_y h$. The tangency condition gives
$$
\langle\operatorname{grad}_x h,\xi\rangle
=-\langle\operatorname{grad}_y h,\eta\rangle.
$$
Taking absolute values and using the fact that $(\xi,\eta)$
has unit length, we obtain
$$
\kappa(x\mid y)\|\xi\|=\kappa(y\mid x)\|\eta\|,
\qquad
\|\xi\|^2+\|\eta\|^2=1.
$$
Solving gives
$$
\|\xi\|
=\frac{\kappa(y\mid x)}
      {\sqrt{\kappa(x\mid y)^2+\kappa(y\mid x)^2}},
\qquad
\|\eta\|
=\frac{\kappa(x\mid y)}
      {\sqrt{\kappa(x\mid y)^2+\kappa(y\mid x)^2}}.
$$

Since $\mathrm d\pi_1(u,v)=u$, the second summand is its kernel.
An orthonormal basis of the first summand together with
$(\xi,\eta)$ maps to mutually orthogonal vectors of lengths
$1,\ldots,1,\|\xi\|$. Their product gives
$$
\mathcal J_d\pi_1(x,y)=\|\xi\|=\frac{\kappa(y\mid x)}
      {\sqrt{\kappa(x\mid y)^2+\kappa(y\mid x)^2}}.
$$
Similarly, 
$$
\mathcal J_d\pi_2(x,y)=\|\eta\|=\frac{\kappa(x\mid y)}
      {\sqrt{\kappa(x\mid y)^2+\kappa(y\mid x)^2}}.
$$
This gives the stated ratio and proves Lemma~\ref{lem:weighted-incidence-jacobian}.
\end{proof}
\begin{lemma}
\label{lem:weighted-incidence-coarea}
For every nonnegative Borel function $q:\Sigma_H\to[0,\infty)$
we have
$$
\int_\M\int_{H_y}q(x,y)
\,\mathrm d\mathcal H^{d-1}(x)\,\mathrm d\operatorname{vol}(y)=
\int_\M\int_{H_x}q(x,y)
\frac{\kappa(x\mid y)}{\kappa(y\mid x)}
\,\mathrm d\mathcal H^{d-1}(y)\,\mathrm d\operatorname{vol}(x).
$$
\end{lemma}

\begin{proof}
By Lemma~\ref{lem:weighted-incidence-structure},
$\Sigma_H$ is a smooth manifold of dimension $2d-1$ and both
coordinate projections $\pi_1,\pi_2:\Sigma_H\to\M$ are smooth
submersions. Since $y\in H_x$ if and only if $x\in H_y$,
their fibers are
$$
\pi_1^{-1}(x)=\{x\}\times H_x,
\qquad
\pi_2^{-1}(y)=H_y\times\{y\}.
$$
Under the product metric the induced volume measures on these
fibers agree with $\mathcal H^{d-1}$ on $H_x$ and $H_y$,
respectively.

Applying the coarea formula \cite[Theorem~2.1]{NicolaescuCoarea2026}
to $\pi_2$ gives
$$
\int_\M\int_{H_y}q(x,y)
\,\mathrm d\mathcal H^{d-1}(x)\,\mathrm d\operatorname{vol}(y)
=\int_{\Sigma_H}q(x,y)\mathcal J_d\pi_2(x,y)
\,\mathrm d\mathcal H^{2d-1}(x,y),
$$
where $\mathcal H^{2d-1}$ denotes the induced volume measure
on $\Sigma_H$. By Lemma~\ref{lem:weighted-incidence-jacobian},
$$
\mathcal J_d\pi_2(x,y)
=\frac{\kappa(x\mid y)}{\kappa(y\mid x)}\mathcal J_d\pi_1(x,y).
$$
After this substitution the coarea formula for $\pi_1$ gives
$$
\begin{aligned}
\int_{\Sigma_H}q(x,y)\mathcal J_d\pi_2(x,y)
\,\mathrm d\mathcal H^{2d-1}(x,y)&=
\int_{\Sigma_H}q(x,y)\frac{\kappa(x\mid y)}{\kappa(y\mid x)}
\mathcal J_d\pi_1(x,y)\,\mathrm d\mathcal H^{2d-1}(x,y)\\
&=
\int_\M\int_{H_x}q(x,y)
\frac{\kappa(x\mid y)}{\kappa(y\mid x)}
\,\mathrm d\mathcal H^{d-1}(y)\,\mathrm d\operatorname{vol}(x).
\end{aligned}
$$
This proves Lemma~\ref{lem:weighted-incidence-coarea}.
\end{proof}
The next lemma gives a uniform bound for the factor
$\kappa(x\mid y)/\kappa(y\mid x)$ in
Lemma~\ref{lem:weighted-incidence-coarea}.

\begin{lemma}
\label{lem:weighted-angle-comparison}
There is a constant $C_*\ge1$ depending only on $(\M,g)$ such
that
\begin{equation}
\label{eq:weighted-angle-comparison}
C_*^{-1}\kappa(y\mid x)
\le \kappa(x\mid y)
\le C_*\kappa(y\mid x)
\end{equation}
for every $(x,y)\in\Sigma_H$.
\end{lemma}

\begin{proof}
Consider all triples $(x,v_1,v_2)$ satisfying
$$
v_1,v_2\in\T_x\M,\qquad
\|v_1\|,\|v_2\|\le D+1,\qquad
\operatorname{Exp}_x(v_1)=\operatorname{Exp}_x(v_2).
$$
These triples form a compact set since the closed disk
bundle of radius $D+1$ is compact and equality of the
endpoints is a closed condition. Let $q_i$ denote the terminal velocity of the geodesic $t\mapsto\operatorname{Exp}_x(tv_i)$ on $[0,1]$ for $i=1,2$.

The ratio $\|q_1-q_2\|/\|v_1-v_2\|$ is continuous wherever
$v_1\ne v_2$. Fix a triple with $v_1=v_2$ and choose bundle
coordinates near its initial and terminal states using
smooth local orthonormal frames. The time-one geodesic flow
is smooth \cite{Paternain1999} and hence Lipschitz on a
sufficiently small coordinate ball.

For each nearby triple, the two initial states have the
same base point, so their Euclidean coordinate distance
is $\|v_1-v_2\|$. Likewise, the two terminal states have
coordinate distance $\|q_1-q_2\|$ because they share the
same endpoint. The Lipschitz estimate therefore bounds
the ratio near each triple with $v_1=v_2$. Compactness gives a constant $C_*\ge1$ depending only on $(\M,g)$ such that $\|q_1-q_2\|\le C_*\|v_1-v_2\|$.
This also holds when $v_1=v_2$, since then $q_1=q_2$.

The reversed geodesics have initial velocities $-q_i$
and terminal velocities $-v_i$. Their speeds are unchanged,
so the same estimate gives
$$
\|v_1-v_2\|
=\|(-v_1)-(-v_2)\|
\le C_*\|(-q_1)-(-q_2)\|
=C_*\|q_1-q_2\|.
$$
Combining the two estimates yields
\begin{equation}
\label{eq:weighted-endpoint-velocity-comparison}
C_*^{-1}\|q_1-q_2\|
\le \|v_1-v_2\|
\le C_*\|q_1-q_2\|.
\end{equation}

For $(x,y)\in\Sigma_H$, take $v_i,q_i$ to be the initial
and terminal velocities of its two minimizing geodesics
parametrized at constant speed on $[0,1]$.
The first variation formula and the definition of $\kappa$
give
$$
\|v_1-v_2\|=d(x,y)\kappa(x\mid y),
\qquad
\|q_1-q_2\|=d(x,y)\kappa(y\mid x).
$$
The conclusion follows from
\eqref{eq:weighted-endpoint-velocity-comparison}.
\end{proof}
\subsection{Averaging and continuity of the singular part}
\label{subsec:weighted-singular-averaging}
\begin{lemma}
\label{lem:weighted-averaged-singular}
For every $x\in U$ and every nonnegative Borel function
$\eta:\M\to[0,\infty]$ we have
\begin{equation}
\label{eq:weighted-J-domination}
\int_{H_x}\eta(y)|w(y)|\psi(y)d(x,y)
\frac{\kappa(x\mid y)^2}{\kappa(y\mid x)}
\,\mathrm d\mathcal H^{d-1}(y)\le C_*^2\|w\|_\infty M_\psi
\int_\M\eta\,\mathrm d\lambda_x.
\end{equation}
In particular, the integral defining $J_{\mu,w}(x)$
converges absolutely for every $x\in U$.
Moreover, on $U$ we have the identity of operator-valued
measures
\begin{equation}
\label{eq:weighted-singular-average}
\int_\M w(y)\,\mathsf S_y|_U\,\mathrm d\mu(y)
=J_{\mu,w}\,\operatorname{vol}|_U.
\end{equation}
\end{lemma}

\begin{proof}
Fix $x\in U$. Since $H_x\subset W_0$, we have
$\psi\le M_\psi$ on $H_x$.
Using \eqref{eq:weighted-angle-comparison} and
Lemma~\ref{lem:weighted-defect}, we obtain
$$
\begin{aligned}
&\int_{H_x}\eta(y)|w(y)|\psi(y)d(x,y)
\frac{\kappa(x\mid y)^2}{\kappa(y\mid x)}
\,\mathrm d\mathcal H^{d-1}(y)\\
&\quad\le
C_*^2\|w\|_\infty M_\psi
\int_{H_x}\eta(y)d(x,y)\kappa(y\mid x)
\,\mathrm d\mathcal H^{d-1}(y)\\
&\quad\le
C_*^2\|w\|_\infty M_\psi
\int_\M\eta\,\mathrm d\lambda_x.
\end{aligned}
$$
This proves \eqref{eq:weighted-J-domination}.
Since $\lambda_x(\M)=K_0$ and
$\|n(x\mid y)\otimes n(x\mid y)\|_{\mathrm{op}}=1$ on $H_x$,
taking $\eta\equiv1$ also proves absolute convergence of
the integral defining $J_{\mu,w}(x)$ with the uniform
absolute integral bound $C_*^2\|w\|_\infty M_\psi K_0$.

For either projection on $\Sigma_H\cap(U\times\M)$,
the fiber integral of a nonnegative Borel function is Borel.
If the integrand is continuous and compactly supported,
its fiber integral is continuous. A countable locally finite partition
of unity with compact supports in submersion charts reduces
these assertions to parameter integrals over $\R^{d-1}$
with smooth volume densities. The first follows from
measurability of parameter integrals. For the second,
only finitely many terms contribute, and each is continuous
by dominated convergence.

In a smooth orthonormal frame on $U$, let $P^{ij}(x,y)$
denote the entries of the smooth operator field
$n(x\mid y)\otimes n(x\mid y)$ on
$\Sigma_H\cap(U\times\M)$.
For any Borel set $E\subset U$, \eqref{eq:weighted-Sy} gives
$$
\mathsf S_y^{ij}(E)=\int_{H_y\cap E}
d(x,y)\kappa(x\mid y)P^{ij}(x,y)
\,\mathrm d\mathcal H^{d-1}(x).
$$
Applying the Borel assertion to the positive and negative
parts on the fibers of $\pi_2$ shows that
$y\mapsto\mathsf S_y^{ij}(E)$ is Borel. Since $|P^{ij}|\le1$,
Lemma~\ref{lem:weighted-defect} gives
$$
|\mathsf S_y^{ij}|(U)\le|\mathsf S_y|(U)\le K_0.
$$
Thus the components of $w(y)\mathsf S_y|_U$ are integrable
signed measure kernels and their integrals define the
left-hand side of \eqref{eq:weighted-singular-average}.

If $x\in U\cap H_y$, then $y\in H_x\subset W_0$.
Since $\mathsf S_y|_U=0$ for $y\notin W_0$, we may restrict
the integration to $y\in W_0$ and replace $\mathrm d\mu(y)$
by $\psi(y)\,\mathrm d\operatorname{vol}(y)$. For each component, apply Lemma~\ref{lem:weighted-incidence-coarea} to the positive and negative parts of the resulting integrand, extended by zero outside $E\times W_0$. This gives
$$
\begin{aligned}
\int_\M w(y)\mathsf S_y^{ij}(E)\,\mathrm d\mu(y)&=\int_E\int_{H_x}w(y)\psi(y)d(x,y)
\frac{\kappa(x\mid y)^2}{\kappa(y\mid x)}
P^{ij}(x,y)\,\mathrm d\mathcal H^{d-1}(y)
\,\mathrm d\operatorname{vol}(x)\\
&=\int_E (J_{\mu,w}(x))^{ij}\,\mathrm d\operatorname{vol}(x).
\end{aligned}
$$
The last equality follows from \eqref{eq:weighted-J-definition}.
This proves \eqref{eq:weighted-singular-average}.
\end{proof}

\begin{lemma}
\label{lem:weighted-J-continuity}
For every $x_0\in U$,
$$
\lim_{x\to x_0}
\bigl\|
\Pi_{x,x_0}J_{\mu,w}(x)\Pi_{x_0,x}
-J_{\mu,w}(x_0)
\bigr\|_{\mathrm{op}}=0.
$$
\end{lemma}

\begin{proof}
Fix $x_0\in U$ and $\varepsilon>0$. The set $Q_{x_0}$ is
compact by Lemma~\ref{lem:weighted-incidence-structure},
and $\lambda_{x_0}(Q_{x_0})=0$ by
Lemma~\ref{lem:weighted-defect}. By outer regularity, choose
a continuous function $0\le\beta\le1$ that equals $1$ on a
neighborhood of $Q_{x_0}$ and satisfies
$\int_\M\beta\,\mathrm d\lambda_{x_0}<\varepsilon$. Then 
$$
J_{\mu,w}=J_{\mu,\beta w}+J_{\mu,(1-\beta) w}.
$$
\eqref{eq:weighted-J-domination} gives $\|J_{\mu,\beta w}(x)\|_{\mathrm{op}}
\le C_*^2\|w\|_\infty M_\psi\int_\M\beta\,\mathrm d\lambda_x$.
Since $\lambda_x\Rightarrow\lambda_{x_0}$ as $x\to x_0$ by
Lemma~\ref{lem:weighted-defect-continuity}, we obtain
\begin{equation}
\label{eq:weighted-cutoff-tail}
\limsup_{x\to x_0}\|J_{\mu,\beta w}(x)\|_{\mathrm{op}}
\le C_*^2\|w\|_\infty M_\psi\int_\M\beta\,\mathrm d\lambda_{x_0}
\le C_*^2\|w\|_\infty M_\psi\varepsilon.
\end{equation}

It remains to show that $J_{\mu,(1-\beta) w}$ is continuous at
$x_0$. Since $\Sigma_H$ is relatively open in the compact
set $\mathcal C$, the set
$$
(\mathcal C\setminus\Sigma_H)
\cap(\M\times\operatorname{supp}(1-\beta))
$$
is compact. Its projection onto the first factor omits $x_0$
because $(1-\beta)$ vanishes on a neighborhood of $Q_{x_0}$.
Choose $\alpha\in C_c^\infty(U)$ equal to $1$ near $x_0$
with support disjoint from this projection. Then
$$
\mathcal C\cap
(\operatorname{supp}\alpha\times\operatorname{supp}(1-\beta))
$$
is a compact subset of $\Sigma_H\cap(U\times\M)$.

In a smooth orthonormal frame on $U$, each entry of the
integrand in \eqref{eq:weighted-J-definition} multiplied by
$\alpha(x)(1-\beta(y))$ is therefore continuous and compactly
supported on $\Sigma_H\cap(U\times\M)$. Indeed, $w$ and
$\psi$ are continuous there as functions of $y$, the normal
projection is smooth, and $\kappa(y\mid x)>0$.
The continuity of compactly supported fiber integrals proved
in Lemma~\ref{lem:weighted-averaged-singular} applies to
$\pi_1$. It follows that $\alpha J_{\mu,(1-\beta) w}$ is continuous
on $U$, hence $J_{\mu,(1-\beta) w}$ is continuous near $x_0$.

The decomposition $J_{\mu,w}=J_{\mu,\beta w}+J_{\mu,(1-\beta)w}$ and the
triangle inequality give
$$
\begin{aligned}
&\bigl\|\Pi_{x,x_0}J_{\mu,w}(x)\Pi_{x_0,x}
      -J_{\mu,w}(x_0)\bigr\|_{\mathrm{op}}\\
&\quad\le \|J_{\mu,\beta w}(x)\|_{\mathrm{op}}+
\bigl\|\Pi_{x,x_0}J_{\mu,(1-\beta)w}(x)\Pi_{x_0,x}
      -J_{\mu,(1-\beta)w}(x_0)\bigr\|_{\mathrm{op}}+\|J_{\mu,\beta w}(x_0)\|_{\mathrm{op}}.
\end{aligned}
$$
The middle term tends to zero by the continuity of
$J_{\mu,(1-\beta)w}$ and the smooth dependence of parallel
transport on $x$.
Using \eqref{eq:weighted-cutoff-tail},
\eqref{eq:weighted-J-domination}, and the choice of $\beta$,
we obtain
$$
\limsup_{x\to x_0}
\bigl\|\Pi_{x,x_0}J_{\mu,w}(x)\Pi_{x_0,x}
      -J_{\mu,w}(x_0)\bigr\|_{\mathrm{op}}\le C_*^2\|w\|_\infty M_\psi
\left(\varepsilon+\int_\M\beta\,\mathrm d\lambda_{x_0}\right)\le 2C_*^2\|w\|_\infty M_\psi\varepsilon.
$$
Letting $\varepsilon\downarrow0$ proves the claim.
\end{proof}

\section{Proof of the volume theorem for moving cut loci}
\label{app:sweep-volume-proof}

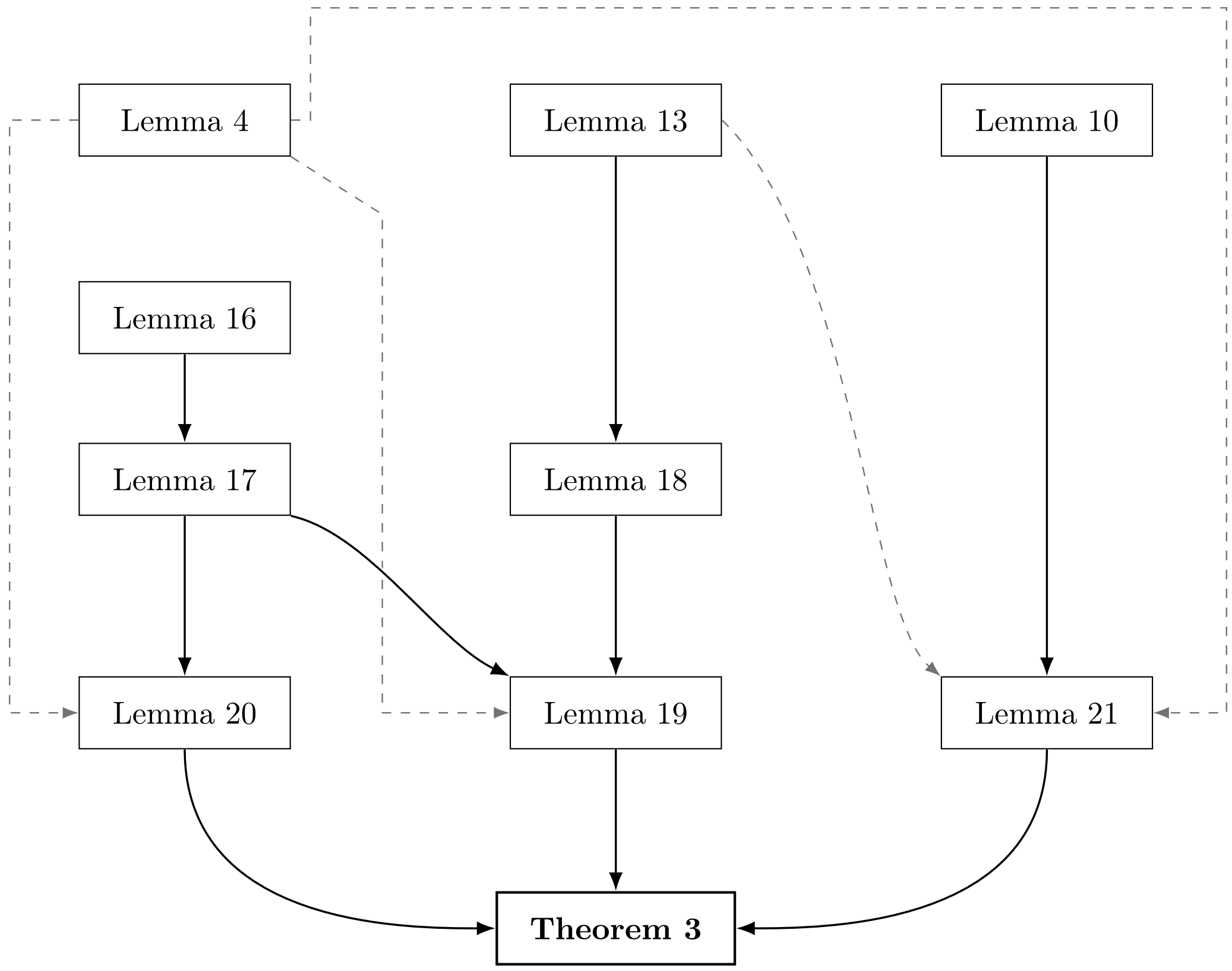
\begin{figure}[H]
\centering
\begin{tikzpicture}[
  >=Latex,
  box/.style={draw=black,rectangle,align=center,
    inner xsep=4pt,inner ysep=5pt,font=\small,
    minimum height=8mm,minimum width=2.35cm},
  main/.style={box,line width=.8pt,font=\small\bfseries,
    minimum width=2.65cm},
  arr/.style={->,line width=.65pt,
    preaction={draw=white,line width=2pt,solid}},
  auxarr/.style={->,dashed,line width=.45pt,draw=black!55}
]
\node[box] (l4) at (-4.8,0) {Lemma~\ref{lem:cut-time-closedness}};
\node[box] (l13) at (0,0) {Lemma~\ref{lem:weighted-angle-comparison}};
\node[box] (l10) at (4.8,0) {Lemma~\ref{lem:weighted-incidence-structure}};
\node[box] (l16) at (-4.8,-2.2) {Lemma~\ref{lem:sweep-upper-support}};
\node[box] (l17) at (-4.8,-4) {Lemma~\ref{lem:sweep-focal-envelope}};
\node[box] (l18) at (0,-4) {Lemma~\ref{lem:sweep-local-transport}};
\node[box] (l20) at (-4.8,-6.6) {Lemma~\ref{lem:sweep-envelope-parametrization}};
\node[box] (l19) at (0,-6.6) {Lemma~\ref{lem:sweep-fixed-envelope}};
\node[box] (l21) at (4.8,-6.6) {Lemma~\ref{lem:sweep-regular-asymptotic}};
\node[main] (t3) at (0,-9) {Theorem~\ref{thm:sweep-volume}};
\draw[auxarr] (l13.east) .. controls (2.65,-1.4) and (2.8,-5.5) .. (l21.north west);
\draw[auxarr] (l4.west) -- (-6.75,0) -- (-6.75,-6.6) -- (l20.west);
\draw[auxarr] (l4.south east) -- (-2.6,-1.05) -- (-2.6,-6.6) -- (l19.west);
\draw[auxarr] (l4.east) -- (-3.4,0) -- (-3.4,1.25) -- (6.8,1.25) -- (6.8,-6.6) -- (l21.east);
\draw[arr] (l16.south) -- (l17.north);
\draw[arr] (l17.south) -- (l20.north);
\draw[arr] (l13.south) -- (l18.north);
\draw[arr] (l18.south) -- (l19.north);
\draw[arr] (l10.south) -- (l21.north);
\draw[arr] (l19.south) -- (t3.north);
\draw[arr] (l17.south east) .. controls (-2.75,-4.6) and (-2,-5.8) .. (l19.north west);
\draw[arr] (l20.south) .. controls (-4.8,-9) and (-2.3,-9) .. (t3.west);
\draw[arr] (l21.south) .. controls (4.8,-9) and (2.3,-9) .. (t3.east);
\end{tikzpicture}
\caption{Proof dependencies for
Theorem~\ref{thm:sweep-volume}. Arrow styles follow Figure~\ref{fig:proof-dependency}.}
\label{fig:sweep-proof-dependency}
\end{figure}

Figure~\ref{fig:sweep-proof-dependency} summarizes the
dependencies among the lemmas in this appendix and
Theorem~\ref{thm:sweep-volume}.

We use the geometric notation of Section~\ref{subsec:notation}.
Throughout the proof we assume $d\ge2$ and treat $d=1$
in the final paragraph.
Set $R=D+1$ and write
$\operatorname{inj}(\M)=\inf_{x\in\M}\operatorname{inj}(x)>0$.
For $E\subset\M$ and $r>0$, write
$$
B_r(E)=\{y\in\M:d(y,E)<r\}
=\bigcup_{z\in E}B_r(z).
$$
In particular, $B_r(\{m\})=B_r(m)$.
We retain $S\M$, $b$ and $\mathcal C$ from
Lemma~\ref{lem:cut-time-closedness} and write
$S_x\M=\{v\in\T_x\M:\|v\|=1\}$.
Let $d_S$ be the intrinsic distance on $S\M$ induced by
the Sasaki metric \cite{Sasaki1958}.
For $(x,v)\in S\M$, let $t_f(x,v)$ be the first conjugate
time along $t\mapsto\operatorname{Exp}_x(tv)$, with value
$\infty$ if there is no conjugate point.
Let $\mathcal C_{\mathrm{conj}}\subset\mathcal C$ consist
of the pairs joined by a minimizing geodesic along which
the endpoints are conjugate.

\begin{lemma}
\label{lem:sweep-upper-support}
Let $\mathcal N$ be a connected Riemannian manifold
with intrinsic distance $d_N$. Let $f\in C(\mathcal N)$
and $K\ge0$. Suppose that every $p\in\mathcal N$ has a neighborhood $V_p$
and a function $s_p:V_p\to\R$ satisfying
$$
s_p(p)=f(p),\qquad
f(q)\le s_p(q),\qquad
|s_p(q)-s_p(p)|\le Kd_N(p,q)
$$
for every $q\in V_p$. Then
$$
|f(q)-f(p)|\le Kd_N(p,q),
\qquad p,q\in\mathcal N.
$$
\end{lemma}

\begin{proof}
Let $\gamma:[0,\ell]\to\mathcal N$ be a piecewise smooth
curve parametrized by arc length.
For each $t<\ell$ we have
$\gamma(t+h)\in V_{\gamma(t)}$ for all sufficiently
small $h>0$. The assumptions therefore give
$$
f(\gamma(t+h))-f(\gamma(t))
\le Kd_N(\gamma(t+h),\gamma(t))
\le Kh.
$$

We claim that
$f(\gamma(\ell))-f(\gamma(0))\le K\ell$.
Otherwise, choose $\varepsilon>0$ such that
$$
f(\gamma(\ell))-f(\gamma(0))>(K+\varepsilon)\ell.
$$
The continuous function
$t\mapsto f(\gamma(t))-(K+\varepsilon)t$
attains its minimum on $[0,\ell]$ at some $t_*<\ell$,
since its value at $\ell$ exceeds its value at $0$.
For all sufficiently small $h>0$, minimality and the
preceding estimate give
$$
(K+\varepsilon)h
\le f(\gamma(t_*+h))-f(\gamma(t_*))
\le Kh.
$$
This is a contradiction and proves the claim.

Applying the same argument to the reversed curve gives
$$
|f(\gamma(\ell))-f(\gamma(0))|\le K\ell.
$$
Taking the infimum over the lengths of piecewise smooth
curves joining $p$ and $q$ proves the result.
\end{proof}
\begin{lemma}
\label{lem:sweep-focal-envelope}
Define
$$
\tau_R(x,v)=\min\{t_f(x,v),R\},
\qquad (x,v)\in S\M,
$$
and
$$
F_m^R=
\left\{\operatorname{Exp}_m\bigl(\tau_R(m,v)v\bigr):
v\in S_m\M\right\}.
$$
Then $\tau_R$ is Lipschitz on $S\M$. Each $F_m^R$ is a
compact Lipschitz image of $S_m\M$ with a Lipschitz constant
independent of $m$.
There is a constant $L_f\ge1$ depending only on $(\M,g)$ such that
$$
(x,y)\in\mathcal C_{\mathrm{conj}},\qquad d(x,m)<\operatorname{inj}(\M)/2
\quad\Longrightarrow\quad
d(y,F_m^R)\le L_f d(x,m).
$$
\end{lemma}

\begin{proof}
We first show that $t_f$ is lower semicontinuous.
Let $(x_j,v_j)\to(x,v)$ and set
$t=\liminf_{j\to\infty}t_f(x_j,v_j)$.
The claim is immediate if $t=\infty$.
Otherwise, after passing to a subsequence, we have
$$
t_f(x_j,v_j)\longrightarrow t
\ge\operatorname{inj}(\M)>0.
$$
Since each
$(\mathrm d\operatorname{Exp}_{x_j})_{t_f(x_j,v_j)v_j}$
is singular, continuity implies that
$(\mathrm d\operatorname{Exp}_x)_{tv}$ is also singular.
Thus $t_f(x,v)\le t$, which proves lower semicontinuity.

By \cite[Theorem~B.3 and Definition~B.1]{FGR2015},
$t_f$ has local upper supports with locally uniform
Lipschitz bounds wherever it is finite.
Thus $\{t_f<\infty\}$ is open and $t_f$ is upper
semicontinuous there.
Together with lower semicontinuity, this gives continuity.
Lemma~\ref{lem:sweep-upper-support} then yields local
Lipschitz regularity on this set.

For $(x,v),(y,w)\in\{t_f<\infty\}$, we have
$$
|\tau_R(x,v)-\tau_R(y,w)|
\le |t_f(x,v)-t_f(y,w)|,
$$
so $\tau_R$ is locally Lipschitz there.
At a point where $t_f=\infty$, lower semicontinuity gives
$t_f>R$ on a neighborhood, where $\tau_R=R$.
Hence $\tau_R$ is locally Lipschitz on all of $S\M$.

Compactness gives a common local Lipschitz constant.
Since $S\M$ is connected for $d\ge2$,
Lemma~\ref{lem:sweep-upper-support} shows that $\tau_R$
is Lipschitz with respect to $d_S$.

The exponential map is smooth on the compact disk bundle
$\{(x,v)\in\T\M:\|v\|\le R\}$.
Together with the Lipschitz continuity of $\tau_R$,
this shows that
$$
(x,v)\longmapsto
\operatorname{Exp}_x\bigl(\tau_R(x,v)v\bigr)
$$
is Lipschitz on $S\M$ with a constant $L_f\ge1$ depending
only on $(\M,g)$.
Its restrictions to $S_m\M$ give the asserted compact
Lipschitz images with a uniform constant.

Now take $(x,y)\in\mathcal C_{\mathrm{conj}}$ and let
$v\in S_x\M$ be the initial direction of a conjugate
minimizing geodesic from $x$ to $y$.
A minimizing geodesic has no conjugate point in its interior.
Since its endpoint is conjugate, we have
$t_f(x,v)=d(x,y)\le D<R$. Hence
$$
y=\operatorname{Exp}_x\bigl(\tau_R(x,v)v\bigr).
$$
If $d(x,m)<\operatorname{inj}(\M)/2$, parallel transport along the short
minimizing geodesic from $x$ to $m$ gives a horizontal curve
in $S\M$ of length $d(x,m)$.
Therefore
$$
d(y,F_m^R)
\le d\left(y,
\operatorname{Exp}_m\bigl(
\tau_R(m,\Pi_{x,m}v)\Pi_{x,m}v\bigr)\right)
\le L_f d_S\bigl((x,v),(m,\Pi_{x,m}v)\bigr)
\le L_f d(x,m).
$$
This proves Lemma~\ref{lem:sweep-focal-envelope}.
\end{proof}

\begin{lemma}
\label{lem:sweep-local-transport}
For every $(x_0,y_0)\in\mathcal C\setminus\mathcal C_{\mathrm{conj}}$,
there are an open product neighborhood $U_0\times W_0$ and
$\eta>0$ such that for every
$(x_a,y_a)\in\mathcal C\cap(U_0\times W_0)$ and every
$x_b\in\M$ with $d(x_a,x_b)<\eta$, there is a point
$y_b\in C_{x_b}$ satisfying
$$
d(y_a,y_b)\le C_*d(x_a,x_b).
$$
Here $C_*$ is the constant in
\eqref{eq:weighted-endpoint-velocity-comparison}.
\end{lemma}

\begin{proof}
The minimizing initial velocities from $x_0$ to $y_0$
form a compact set. Nonconjugacy and the inverse function
theorem imply that each is isolated, so their number $N$
is finite.

Applying the inverse function theorem to
$(x,v)\mapsto(x,\operatorname{Exp}_x(v))$ gives smooth
branches $v_1,\ldots,v_N$. Choose strongly convex
neighborhoods $U$ of $x_0$ and $W$ of $y_0$ so that
the branches are defined on a neighborhood of
$\overline U\times\overline W$, remain pairwise distinct,
have nonconjugate endpoints, and have speeds less than $R$.

We may also require that every minimizing geodesic between
points of $U\times W$ belongs to one of these branches.
Otherwise, there would be point pairs approaching
$(x_0,y_0)$ with minimizing initial conditions outside
all the chosen inverse function neighborhoods.
By compactness, a subsequence would converge to a minimizing
initial condition from $x_0$ to $y_0$ and hence eventually
enter one of those neighborhoods. This is a contradiction.

Set $f_i(x,y)=\frac12\|v_i(x,y)\|^2$. Then
$$
\frac12d(x,y)^2=\min_{1\le i\le N}f_i(x,y)
\qquad\text{on }U\times W.
$$
Since all minimizing branches are nonconjugate, a point
pair belongs to $\mathcal C$ exactly when at least two
branches attain the minimum. The first variation formula
for energy and
\eqref{eq:weighted-endpoint-velocity-comparison} give
$$
\operatorname{grad}_y(f_i-f_j)\ne0,
\qquad
\|\operatorname{grad}_x(f_i-f_j)\|
\le C_*\|\operatorname{grad}_y(f_i-f_j)\|,
\qquad i\ne j.
$$
The terminal velocities at the common endpoint are distinct
by uniqueness of the geodesic flow.

Choose open neighborhoods $U_0\Subset U$ of $x_0$
and $W_0\Subset W$ of $y_0$.
Take
$$
0<\eta<
\frac12\min\left\{
\operatorname{inj}(\M),
d(\overline U_0,\M\setminus U),
\frac{d(\overline W_0,\M\setminus W)}{C_*}
\right\}.
$$
If $y_a\in C_{x_b}$, take $y_b=y_a$.
Otherwise, let $i$ be the unique minimizing branch at
$(x_b,y_a)$ and choose a minimizing branch $j\ne i$
at $(x_a,y_a)$. Then
$$
(f_i-f_j)(x_a,y_a)\ge0,
\qquad
(f_i-f_j)(x_b,y_a)<0.
$$

Let $x(t)$ be the unit speed minimizing geodesic from $x_a$
to $x_b$ on $[0,\ell]$ where $\ell=d(x_a,x_b)<\eta$.
Consider the initial value problem
$$
\dot y(t)
=-\frac{\mathrm d_x(f_i-f_j)(x(t),y(t))[\dot x(t)]}
        {\|\operatorname{grad}_y(f_i-f_j)(x(t),y(t))\|^2}
  \operatorname{grad}_y(f_i-f_j)(x(t),y(t)),
\qquad y(0)=y_a.
$$
The right-hand side is smooth, so a local solution exists.
The Cauchy--Schwarz inequality and the preceding gradient
estimate give
$$
\|\dot y(t)\|
=\frac{
|\mathrm d_x(f_i-f_j)(x(t),y(t))[\dot x(t)]|
}{
\|\operatorname{grad}_y(f_i-f_j)(x(t),y(t))\|
}\le
\frac{
\|\operatorname{grad}_x(f_i-f_j)(x(t),y(t))\|
\,\|\dot x(t)\|
}{
\|\operatorname{grad}_y(f_i-f_j)(x(t),y(t))\|
}
\le C_*,
$$
where we used $\|\dot x(t)\|=1$. We also have
$$
\begin{aligned}
\frac{\mathrm d}{\mathrm dt}(f_i-f_j)(x(t),y(t))
=\mathrm d_x(f_i-f_j)(x(t),y(t))[\dot x(t)]+\left\langle
\operatorname{grad}_y(f_i-f_j)(x(t),y(t)),\dot y(t)
\right\rangle
=0.
\end{aligned}
$$

For as long as the solution exists with $0\le t\le\ell$,
$$
\begin{aligned}
d(x_a,x(t))=t\le\ell<\eta,\qquad
d(y_a,y(t))
\le\int_0^t\|\dot y(s)\|\,\mathrm ds
\le C_*t\le C_*\ell<C_*\eta.
\end{aligned}
$$
These estimates and the choice of $\eta$ keep both curves
in compact subsets of $U$ and $W$, respectively.
The denominator in the differential equation is bounded
away from zero on $\overline U\times\overline W$.
The continuation theorem for ordinary differential equations
therefore extends the solution to $[0,\ell]$.

Set $z=y(\ell)$. Then
$$
d(y_a,z)\le C_*\ell,
\qquad
(f_i-f_j)(x_b,z)=(f_i-f_j)(x_a,y_a)\ge0.
$$
Thus branch $i$ is not the unique minimizing branch at
$(x_b,z)$.
Follow the minimizing geodesic in $W$ from $y_a$ to $z$
and let $y_b$ be the first point where branch $i$ ceases
to be uniquely minimizing.
Such a point exists because branch $i$ is uniquely
minimizing at $(x_b,y_a)$ but not at $(x_b,z)$.
Continuity gives
$$
f_i(x_b,y_b)
=\min_{k\ne i}f_k(x_b,y_b)
=\min_{1\le k\le N}f_k(x_b,y_b).
$$
Hence at least two distinct branches minimize at $(x_b,y_b)$,
so $y_b\in C_{x_b}$.
Since $y_b$ lies on the minimizing geodesic from $y_a$ to $z$,
$$
d(y_a,y_b)\le d(y_a,z)\le C_*\ell.
$$
\end{proof}

\begin{lemma}
\label{lem:sweep-fixed-envelope}
There are constants $C_{\mathrm{env}}>0$ and $\delta_0>0$
depending only on $(\M,g)$ such that
$$
\mathcal A_\delta(m)
\subset B_{C_{\mathrm{env}}\delta}(E_m),
\qquad E_m=C_m\cup F_m^R,
$$
for every $m\in\M$ and $0<\delta<\delta_0$.
\end{lemma}

\begin{proof}
Choose $\delta_0<\operatorname{inj}(\M)/2$ and
$C_{\mathrm{env}}>C_*+L_f+2$. Fix $x\in B_\delta(m)$ and
$y\in C_x$. The assertion is immediate if $x=m$ or $d(y,F_m^R)\le(C_*+L_f+2)\delta$.
Assume otherwise. Let $x(t)$ be the unit speed minimizing
geodesic from $x$ to $m$ on $[0,\ell]$, where
$\ell=d(x,m)<\delta$, and consider
$$
Z=\left\{(t,z):0\le t\le\ell,\ z\in C_{x(t)},
                 \ d(z,y)\le(C_*+1)\delta\right\}.
$$
This set is compact by the closedness of $\mathcal C$ in
Lemma~\ref{lem:cut-time-closedness}. If a corresponding
point pair $(x(t),z)$ belonged to $\mathcal C_{\mathrm{conj}}$,
Lemma~\ref{lem:sweep-focal-envelope} would give
$$
d(y,F_m^R)\le d(y,z)+d(z,F_m^R)\le(C_*+1)\delta+L_fd(x(t),m)\le(C_*+1+L_f)\delta,
$$
which is a contradiction. Thus every corresponding cut
pair is nonconjugate.

Cover the compact set of these point pairs by finitely
many product neighborhoods from
Lemma~\ref{lem:sweep-local-transport}.
Choose a partition $0=t_0<\cdots<t_N=\ell$
with increments smaller than all the associated
admissible step sizes.
Starting with $y_0=y$, apply that lemma successively
to choose $y_j\in C_{x(t_j)}$ satisfying
$$
d(y_{j-1},y_j)\le C_*(t_j-t_{j-1}).
$$
At every step,
$$
d(y,y_j)
\le\sum_{k=1}^j d(y_{k-1},y_k)
\le C_*t_j<C_*\delta.
$$
Hence $(t_j,y_j)\in Z$, so the finite cover remains
applicable and the construction continues to $t_N=\ell$.
We obtain $y_N\in C_m$ with
$$
d(y,y_N)\le C_*\ell<C_*\delta,
$$
which proves the inclusion. 
\end{proof}

\begin{lemma}
\label{lem:sweep-envelope-parametrization}
Let $E_m$ be the fixed envelope defined in
Lemma~\ref{lem:sweep-fixed-envelope}. Then the following
statements hold.

\textup{(i)} For every $m\in\M$, both $C_m$ and $F_m^R$ are
Lipschitz images of $S_m\M$ with constants uniform in $m$.
Consequently, $E_m$ is a finite union of Lipschitz images of
compact subsets of $\R^{d-1}$ and $\mathcal H^{d-1}(E_m)<\infty$.
The number of parameter domains can be chosen independently
of $m$.

\textup{(ii)} There exist constants $C_1,r_0>0$ depending only
on $(\M,g)$ such that
\begin{equation}
\label{eq:sweep-envelope-tube-bound}
\operatorname{vol}(B_r(E_m))\le C_1r,
\qquad m\in\M,\quad 0<r<r_0.
\end{equation}

\textup{(iii)} There exists a constant $C_2>0$ depending only
on $(\M,g)$ such that for every $m\in\M$ and every compact
set $K\subset E_m$ we have
$$
\limsup_{r\downarrow0}
\frac{\operatorname{vol}(B_r(K))}{r}
\le C_2\mathcal H^{d-1}(K).
$$
\end{lemma}

\begin{proof}
By \cite[Theorem~B.2(1)]{FGR2015}, the function
$b(m,\cdot)$ has local upper supports with a Lipschitz
bound independent of $m$.
Since $b$ is continuous, Lemma~\ref{lem:sweep-upper-support}
and comparison of great-circle arc length with chord length
give
$$
|b(m,v)-b(m,w)|
\le C\arccos\langle v,w\rangle
\le\frac{\pi C}{2}\|v-w\|,
\qquad v,w\in S_m\M,
$$
where $C$ is independent of $m$.

Since $b\le D$ and $\|w\|=1$,
$$
\begin{aligned}
\|b(m,v)v-b(m,w)w\|
&=\left\|
b(m,v)(v-w)+\bigl(b(m,v)-b(m,w)\bigr)w
\right\|\\
&\le D\|v-w\|+|b(m,v)-b(m,w)|\\
&\le\left(D+\frac{\pi C}{2}\right)\|v-w\|.
\end{aligned}
$$
The exponential map has uniformly bounded derivative
on the compact disk bundle of radius $D$.
Consequently,
$$
v\longmapsto\operatorname{Exp}_m(b(m,v)v)
$$
is a Lipschitz map onto $C_m$ with a bound independent of $m$.

For $v,w\in S_m\M$, the great-circle arc joining them
has the same length in the Sasaki metric, so
$$
d_S\bigl((m,v),(m,w)\bigr)
\le\arccos\langle v,w\rangle
\le\frac{\pi}{2}\|v-w\|.
$$
Together with Lemma~\ref{lem:sweep-focal-envelope},
this shows that
$$
v\longmapsto
\operatorname{Exp}_m\bigl(\tau_R(m,v)v\bigr)
$$
is a Lipschitz map onto $F_m^R$ with a bound independent of $m$.

Covering the unit sphere by a fixed number of compact
coordinate patches gives
$$
E_m=\bigcup_{\alpha=1}^{A}T_\alpha(P_\alpha),
$$
where each $P_\alpha\subset\R^{d-1}$ is compact and each
$T_\alpha:P_\alpha\to\M$ is Lipschitz. The number $A$ and
the Lipschitz bounds can be chosen independently of $m$.
Since each $P_\alpha\subset\R^{d-1}$ is compact,
the Hausdorff measure estimate for Lipschitz maps gives
$$
\mathcal H^{d-1}(E_m)
\le\sum_{\alpha=1}^{A}
\mathcal H^{d-1}\bigl(T_\alpha(P_\alpha)\bigr)
<\infty.
$$
This proves \textup{(i)}.

The unit sphere can be covered by $O(r^{-(d-1)})$ balls
of radius $r$. The uniform Lipschitz bounds for the two
parametrizations therefore give a cover of $E_m$ by
$O(r^{-(d-1)})$ balls of radius $O(r)$ in $\M$.
Enlarging their radii by $r$ gives a cover of $B_r(E_m)$.
Each enlarged ball has volume $O(r^d)$ for sufficiently
small $r$, uniformly in its center. Summing over the cover
gives
$$
\operatorname{vol}(B_r(E_m))
\le C r^{-(d-1)}r^d
=Cr,
\qquad m\in\M,\quad 0<r<r_0,
$$
where $C$ and $r_0$ depend only on $(\M,g)$.
This proves \textup{(ii)} with $C_1=C$.

To prove \textup{(iii)}, we use the Minkowski content theorem
\cite[Theorem~4.29]{Mattila2023}. It states that every compact
set $F\subset\R^d$ which is the Lipschitz image of a bounded
subset of $\R^{d-1}$ satisfies
$$
\lim_{r\downarrow0}
\frac{\mathcal L^d(B_r(F))}{2r}
=\mathcal H^{d-1}(F).
$$
Here $\mathcal L^d$ denotes Lebesgue measure on $\R^d$,
and the neighborhoods and Hausdorff measure are taken
with respect to the Euclidean metric.

Choose finitely many coordinate neighborhoods $V_j$
with coordinate maps $\xi_j:V_j\to\R^d$ and compact
subsets $K_j^0\Subset V_j$ covering $\M$.
On a fixed compact neighborhood of each $K_j^0$
contained in $V_j$, choose $c_j\ge1$ such that
$$
\|\xi_j(x)-\xi_j(y)\|\le c_jd(x,y),
\qquad
\operatorname{vol}(E)
\le c_j\mathcal L^d(\xi_j(E))
$$
for all points $x,y$ and Borel sets $E$ in that neighborhood.
Such constants exist by smoothness of the coordinate maps
and boundedness of the volume densities.
They depend only on the fixed charts and the metric.

For a compact set $K\subset E_m$, put
$$
K_{\alpha j}
=K\cap T_\alpha(P_\alpha)\cap K_j^0.
$$
Each $K_{\alpha j}$ is compact.
The set $T_\alpha^{-1}(K_{\alpha j})$ is a compact
subset of $P_\alpha\subset\R^{d-1}$, and
$$
\xi_j(K_{\alpha j})
=(\xi_j\circ T_\alpha)
\bigl(T_\alpha^{-1}(K_{\alpha j})\bigr).
$$
The composition is Lipschitz on this preimage.
Hence the Minkowski content theorem applies to
$\xi_j(K_{\alpha j})$.

For sufficiently small $r$, the neighborhood
$B_r(K_{\alpha j})$ stays in the fixed neighborhood
on which the coordinate estimates hold.
The distance comparison gives
$$
\xi_j(B_r(K_{\alpha j}))
\subset B_{c_jr}(\xi_j(K_{\alpha j})).
$$
Consequently,
$$
\begin{aligned}
\operatorname{vol}(B_r(K_{\alpha j}))
\le c_j\mathcal L^d
\bigl(\xi_j(B_r(K_{\alpha j}))\bigr)
\le c_j\mathcal L^d
\bigl(B_{c_jr}(\xi_j(K_{\alpha j}))\bigr).
\end{aligned}
$$
Applying the Minkowski content theorem and the Hausdorff
measure estimate for the Lipschitz map $\xi_j$ yields
$$
\begin{aligned}
\limsup_{r\downarrow0}
\frac{\operatorname{vol}(B_r(K_{\alpha j}))}{r}
\le c_j\lim_{r\downarrow0}
\frac{
\mathcal L^d(B_{c_jr}(\xi_j(K_{\alpha j})))
}{r}
=2c_j^2\mathcal H^{d-1}(\xi_j(K_{\alpha j}))
\le2c_j^{d+1}\mathcal H^{d-1}(K_{\alpha j}).
\end{aligned}
$$

Since $K=\bigcup_{\alpha,j}K_{\alpha j}$, we have
$$
B_r(K)=\bigcup_{\alpha,j}B_r(K_{\alpha j}).
$$
Finite subadditivity therefore gives
$$
\begin{aligned}
\limsup_{r\downarrow0}
\frac{\operatorname{vol}(B_r(K))}{r}
\le
2\sum_{\alpha,j}c_j^{d+1}
\mathcal H^{d-1}(K_{\alpha j})
\le
2A\left(\sum_jc_j^{d+1}\right)
\mathcal H^{d-1}(K).
\end{aligned}
$$
The number $A$ and the constants $c_j$ are independent
of $m$ and $K$.
This proves \textup{(iii)} with
$C_2=2A\sum_jc_j^{d+1}$.
\end{proof}

\begin{lemma}
\label{lem:sweep-regular-asymptotic}
Fix $m\in\M$. If $f\in C(\M)$ satisfies
$\operatorname{supp}f\cap C_m\subset H_m$, then
$$
\lim_{\delta\downarrow0}\frac1\delta
\int_{\mathcal A_\delta(m)}f(y)\,\mathrm d\operatorname{vol}(y)
=2\int_{H_m}f(y)\frac{\kappa(m\mid y)}{\kappa(y\mid m)}\,\mathrm d\mathcal H^{d-1}(y).
$$
\end{lemma}

\begin{proof}
If $\operatorname{supp}f\cap C_m=\varnothing$,
Lemma~\ref{lem:cut-time-closedness}\textup{(iii)} gives
$$
\mathcal A_\delta(m)
\subset\M\setminus\operatorname{supp}f
$$
for all sufficiently small $\delta$.
Both integrals then vanish because $H_m\subset C_m$.

Fix $y_0\in H_m$.
By Lemma~\ref{lem:weighted-incidence-structure},
there are $0<\rho<\operatorname{inj}(\M)$,
an open neighborhood $W$ of $y_0$, and a smooth difference
$h$ of two geodesic length branches such that
\begin{equation}
\label{eq:sweep-local-cut-equivalence}
y\in C_x\Longleftrightarrow h(x,y)=0,
\qquad (x,y)\in B_\rho(m)\times W.
\end{equation}
For $(x,y)\in B_\rho(m)\times W$ with $h(x,y)=0$,
$$
\|\operatorname{grad}_x h(x,y)\|=\kappa(x\mid y),
\qquad
\|\operatorname{grad}_y h(x,y)\|=\kappa(y\mid x).
$$

Choose a relatively compact open patch $S$ of $H_m\cap W$
containing $y_0$ whose closure lies in a larger regular patch.
Fix the order of the two branches on that larger patch so that
$$
n(z\mid m)=\frac{\operatorname{grad}_y h(m,z)}{\kappa(z\mid m)}
$$
is a smooth unit normal. We may choose $S$ and $a>0$ so that
$$
\Theta:S\times(-a,a)\longrightarrow\mathcal T\subset W,
\qquad
\Theta(z,t)=\operatorname{Exp}_z(tn(z\mid m))
$$
is a diffeomorphism and extends smoothly to a neighborhood of
$\overline S\times[-a,a]$. Since
$$
\partial_t h(m,\Theta(z,0))=\kappa(z\mid m)>0,
$$
compactness of $\overline S$ allows us to take $a$ and
$\rho$ sufficiently small that
$$
\partial_t h(\operatorname{Exp}_m(u),\Theta(z,t))>0
$$
for all $u\in\T_m\M$, $z\in\overline S$ and $t\in\R$
satisfying $\|u\|<\rho$ and $|t|\le a$.
The implicit function theorem and this strict monotonicity
give a smooth function $\phi(u,z)$ such that
\begin{equation}
\label{eq:sweep-unique-normal-graph}
\phi(0,z)=0,\qquad |\phi(u,z)|<a/2,\qquad
h(\operatorname{Exp}_m(u),\Theta(z,t))=0
\Longleftrightarrow t=\phi(u,z)
\end{equation}
whenever $\|u\|<\rho$, $z\in\overline S$ and $|t|<a$.
Here $\rho$ can be chosen uniformly for $z\in\overline S$:
local implicit functions exist on a finite cover of
$\overline S$ and agree wherever their domains overlap by
strict monotonicity.

For $0<\delta<\rho/2$, set
$$
I_\delta(z)=\{\phi(u,z):u\in\T_m\M,\ \|u\|<\delta\}.
$$
Then
\begin{equation}
\label{eq:sweep-exact-tube-set}
\mathcal A_\delta(m)\cap\mathcal T
=\{\Theta(z,t):z\in S,\ t\in I_\delta(z)\}.
\end{equation}
Indeed, $\Theta(z,t)$ belongs to the left-hand side exactly
when it lies in $C_x$ for some $x\in B_\delta(m)$.
Such an $x$ has the unique representation
$x=\operatorname{Exp}_m(u)$ with $\|u\|<\delta$.
Equations~\eqref{eq:sweep-local-cut-equivalence} and
\eqref{eq:sweep-unique-normal-graph} then make this condition
equivalent to $t\in I_\delta(z)$.

Differentiating $h(\operatorname{Exp}_m(u),\Theta(z,\phi(u,z)))=0$
at $u=0$ in the direction $V\in\T_m\M$ gives
$$
\langle\operatorname{grad}_x h(m,z),V\rangle
+\kappa(z\mid m)\,\mathrm d_u\phi(0,z)V=0,
$$
where we used $(\mathrm d\operatorname{Exp}_m)_0=\operatorname{Id}$
and $\partial_t h(m,\Theta(z,0))=\kappa(z\mid m)$.
Hence
$$
\mathrm d_u\phi(0,z)V
=-\frac{\langle\operatorname{grad}_x h(m,z),V\rangle}
        {\kappa(z\mid m)}.
$$
The second derivatives of $\phi$ with respect to $u$
are uniformly bounded on the compact set
$\{u\in\T_m\M:\|u\|\le\rho/2\}\times\overline S$.
Taylor's theorem therefore gives
$$
\phi(u,z)
=-\frac{\langle\operatorname{grad}_x h(m,z),u\rangle}
        {\kappa(z\mid m)}
+O_S(\|u\|^2),
\qquad z\in\overline S,\quad \|u\|\le\rho/2.
$$
For each $z$, the continuous image $I_\delta(z)$ of the
connected open ball is an interval.
Continuity of $\phi(\cdot,z)$ and density of the open ball
in the closed ball give
$$
\sup I_\delta(z)=\max_{\|u\|\le\delta}\phi(u,z),
\qquad
\inf I_\delta(z)=\min_{\|u\|\le\delta}\phi(u,z).
$$
The Cauchy--Schwarz inequality gives
$$
\max_{\|u\|\le\delta}
\left(
-\frac{\langle\operatorname{grad}_x h(m,z),u\rangle}
       {\kappa(z\mid m)}
\right)
=\delta\frac{\kappa(m\mid z)}{\kappa(z\mid m)},
$$
with equality when $u$ points opposite to
$\operatorname{grad}_x h(m,z)$ and has length $\delta$.
The minimum is the negative of this value.
The uniform Taylor remainder changes each extremum
by at most $O_S(\delta^2)$, so
$$
\sup I_\delta(z)=\delta\frac{\kappa(m\mid z)}{\kappa(z\mid m)}+O_S(\delta^2),
\qquad
\inf I_\delta(z)=-\delta\frac{\kappa(m\mid z)}{\kappa(z\mid m)}+O_S(\delta^2).
$$
Consequently, for all sufficiently small $\delta$,
\begin{equation}
\label{eq:sweep-uniform-interval-width}
\begin{aligned}
|I_\delta(z)|
=2\delta\frac{\kappa(m\mid z)}{\kappa(z\mid m)}
+O_S(\delta^2),\qquad I_\delta(z)\subset[-2C_*\delta,2C_*\delta].
\end{aligned}
\end{equation}
By Lemma~\ref{lem:weighted-angle-comparison},
$$
\frac{\kappa(m\mid z)}{\kappa(z\mid m)}\le C_*.
$$
The uniform Taylor estimate therefore gives a constant
$C$ independent of $z$ such that
$$
\sup_{\|u\|<\delta}|\phi(u,z)|
\le C_*\delta+C\delta^2
\le2C_*\delta
$$
for all sufficiently small $\delta$.
This proves the inclusion, and both estimates are uniform
for $z\in\overline S$.

Since $\Theta$ is smooth near $\overline S\times[-a,a]$
and its volume Jacobian equals $1$ at $t=0$,
we have the following expansion uniformly for $z\in\overline S$:
$$
\mathrm d\operatorname{vol}(\Theta(z,t))
=\bigl(1+O_S(|t|)\bigr)
\,\mathrm dt\,\mathrm d\mathcal H^{d-1}(z).
$$

Suppose that $f$ has compact support in $\mathcal T$.
By \eqref{eq:sweep-exact-tube-set} and Fubini's theorem,
\begin{equation}
\label{eq:sweep-tube-fubini}
\int_{\mathcal A_\delta(m)}f\,\mathrm d\operatorname{vol}
=\int_S\int_{I_\delta(z)}
f(\Theta(z,t))
\,\mathrm dt\,\mathrm d\mathcal H^{d-1}(z)
+O(\delta^2).
\end{equation}
Indeed, \eqref{eq:sweep-uniform-interval-width} and
$\mathcal H^{d-1}(S)<\infty$ give
$$
\int_S\int_{I_\delta(z)}
|t|\,\mathrm dt\,\mathrm d\mathcal H^{d-1}(z)
\le2C_*\delta
\int_S|I_\delta(z)|\,\mathrm d\mathcal H^{d-1}(z)
=O(\delta^2).
$$
Since $f$ is bounded and the volume density remainder
is uniformly $O_S(|t|)$, this gives the stated error estimate.

Uniform continuity of $f$ gives
$$
\sup_{\substack{z\in\overline S\\ |t|\le2C_*\delta}}
|f(\Theta(z,t))-f(z)|
\longrightarrow0.
$$
Since $\mathcal H^{d-1}(S)<\infty$, equation
\eqref{eq:sweep-uniform-interval-width} implies
$$
\int_{\mathcal A_\delta(m)}f\,\mathrm d\operatorname{vol}
=\int_S f(z)|I_\delta(z)|\,\mathrm d\mathcal H^{d-1}(z)
+o(\delta)
=2\delta\int_S f(z)\frac{\kappa(m\mid z)}{\kappa(z\mid m)}\,\mathrm d\mathcal H^{d-1}(z)
+o(\delta).
$$
As $H_m\cap\mathcal T=S$, this proves
\begin{equation}
\label{eq:sweep-regular-test-limit}
\lim_{\delta\downarrow0}\frac1\delta
\int_{\mathcal A_\delta(m)}f\,\mathrm d\operatorname{vol}
=2\int_{H_m}f(z)\frac{\kappa(m\mid z)}{\kappa(z\mid m)}\,\mathrm d\mathcal H^{d-1}(z).
\end{equation}

For a general $f$ in the statement, the compact set
$\operatorname{supp}f\cap C_m$ is covered by $N$ such
normal tubes.
Multiplying $f$ by a finite smooth partition of unity
on a neighborhood of $\operatorname{supp}f$ subordinate
to these tubes and $\M\setminus C_m$ gives
$$
f=f_0+\sum_{j=1}^{N}f_j,
$$
where $f_j$ has compact support in the $j$th tube
for $1\le j\le N$, and $f_0$ has compact support
disjoint from $C_m$.

By Lemma~\ref{lem:cut-time-closedness}\textup{(iii)},
the integral of $f_0$ over $\mathcal A_\delta(m)$
vanishes for all sufficiently small $\delta$.
Also, $f_0=0$ on $H_m$.
Applying \eqref{eq:sweep-regular-test-limit} to each $f_j$
and summing the finitely many limits yields
$$
\begin{aligned}
\lim_{\delta\downarrow0}\frac1\delta
\int_{\mathcal A_\delta(m)}f\,\mathrm d\operatorname{vol}
&=\sum_{j=1}^{N}
\lim_{\delta\downarrow0}\frac1\delta
\int_{\mathcal A_\delta(m)}f_j\,\mathrm d\operatorname{vol}\\
&=2\sum_{j=1}^{N}\int_{H_m}
f_j(z)\frac{\kappa(m\mid z)}{\kappa(z\mid m)}
\,\mathrm d\mathcal H^{d-1}(z)\\
&=2\int_{H_m}
f(z)\frac{\kappa(m\mid z)}{\kappa(z\mid m)}
\,\mathrm d\mathcal H^{d-1}(z).
\end{aligned}
$$
This proves the lemma.
\end{proof}

\begin{proof}[Proof of Theorem~\ref{thm:sweep-volume}]
First suppose that $d\ge2$.
Lemmas~\ref{lem:sweep-fixed-envelope} and
\ref{lem:sweep-envelope-parametrization} allow us to choose
$\delta_0>0$ depending only on $(\M,g)$ so that
$C_{\mathrm{env}}\delta_0<r_0$ and
$$
\operatorname{vol}(\mathcal A_\delta(m))
\le\operatorname{vol}(B_{C_{\mathrm{env}}\delta}(E_m))
\le C_1C_{\mathrm{env}}\delta,
\qquad m\in\M,\quad 0<\delta<\delta_0.
$$
This proves \textup{(i)}.

For \textup{(ii)}, fix $m\in\M$ and $f\in C(\M)$.
Since $H_m\subset E_m$,
Lemmas~\ref{lem:sweep-envelope-parametrization}
and~\ref{lem:weighted-angle-comparison} give
$$
\int_{H_m}|f(y)|\frac{\kappa(m\mid y)}{\kappa(y\mid m)}\,\mathrm d\mathcal H^{d-1}(y)
\le C_*\|f\|_\infty\mathcal H^{d-1}(E_m)<\infty.
$$
If $Q_m=\varnothing$, the asserted limit follows directly
from Lemma~\ref{lem:sweep-regular-asymptotic}.
Suppose now that $Q_m\ne\varnothing$.
The set $Q_m$ is compact because $H_m$ is relatively open
in $C_m$ by Lemma~\ref{lem:weighted-incidence-structure}.
For $\varepsilon>0$, set
$$
K_\varepsilon
=E_m\cap\{z\in\M:d(z,Q_m)\le3\varepsilon\}.
$$
When $C_{\mathrm{env}}\delta<\varepsilon$,
Lemma~\ref{lem:sweep-fixed-envelope} gives
$$
\mathcal A_\delta(m)\cap B_{2\varepsilon}(Q_m)
\subset B_{C_{\mathrm{env}}\delta}(K_\varepsilon).
$$
Indeed, every point $y$ on the left lies within
$C_{\mathrm{env}}\delta$ of some $z\in E_m$.
Then $d(z,Q_m)<3\varepsilon$, so $z\in K_\varepsilon$.
Applying the compact-set estimate in
Lemma~\ref{lem:sweep-envelope-parametrization} with
$r=C_{\mathrm{env}}\delta$ yields
\begin{equation}
\label{eq:sweep-exceptional-tail-bound}
\limsup_{\delta\downarrow0}
\frac{\operatorname{vol}
(\mathcal A_\delta(m)\cap B_{2\varepsilon}(Q_m))}{\delta}
\le C_{\mathrm{env}}C_2\mathcal H^{d-1}(K_\varepsilon).
\end{equation}
The compact sets $K_\varepsilon$ decrease to $Q_m$ as
$\varepsilon\downarrow0$.
Since $\mathcal H^{d-1}(E_m)<\infty$ and
$\mathcal H^{d-1}(Q_m)=0$, continuity from above gives
$$
\mathcal H^{d-1}(K_\varepsilon)\longrightarrow0.
$$
In particular,
$$
\lim_{\varepsilon\downarrow0}\limsup_{\delta\downarrow0}
\frac{\operatorname{vol}
(\mathcal A_\delta(m)\cap B_{2\varepsilon}(Q_m))}{\delta}=0.
$$

Choose $\chi_\varepsilon\in C^\infty(\M)$ such that
$$
0\le\chi_\varepsilon\le1,\qquad
\chi_\varepsilon|_{B_\varepsilon(Q_m)}=0,\qquad
\chi_\varepsilon|_{\M\setminus B_{2\varepsilon}(Q_m)}=1.
$$
For fixed $\varepsilon>0$, the set
$\operatorname{supp}(f\chi_\varepsilon)\cap C_m$ is a
compact subset of $H_m$.
Lemma~\ref{lem:sweep-regular-asymptotic} therefore gives
$$
\lim_{\delta\downarrow0}\frac1\delta
\int_{\mathcal A_\delta(m)}f\chi_\varepsilon
\,\mathrm d\operatorname{vol}
=2\int_{H_m}f(y)\chi_\varepsilon(y)\frac{\kappa(m\mid y)}{\kappa(y\mid m)}
\,\mathrm d\mathcal H^{d-1}(y).
$$
On the other hand, \eqref{eq:sweep-exceptional-tail-bound}
implies
$$
\limsup_{\delta\downarrow0}\frac1\delta
\left|\int_{\mathcal A_\delta(m)}f(1-\chi_\varepsilon)
\,\mathrm d\operatorname{vol}\right|
\le C_{\mathrm{env}}C_2\|f\|_\infty
\mathcal H^{d-1}(K_\varepsilon).
$$
For every $y\in H_m$, closedness of $Q_m$ implies
$d(y,Q_m)>0$, so $\chi_\varepsilon(y)=1$ for all sufficiently
small $\varepsilon$.
Dominated convergence and the absolute integrability proved
above show that
$$
\int_{H_m}|f(y)|(1-\chi_\varepsilon(y))\frac{\kappa(m\mid y)}{\kappa(y\mid m)}
\,\mathrm d\mathcal H^{d-1}(y)\longrightarrow0.
$$
These estimates also apply when $f$ changes sign. In fact,
$$
\begin{aligned}
&\limsup_{\delta\downarrow0}
\left|\frac1\delta\int_{\mathcal A_\delta(m)}f
\,\mathrm d\operatorname{vol}
-2\int_{H_m}f(y)\frac{\kappa(m\mid y)}{\kappa(y\mid m)}
\,\mathrm d\mathcal H^{d-1}(y)\right|\\
&\quad\le C_{\mathrm{env}}C_2\|f\|_\infty
\mathcal H^{d-1}(K_\varepsilon)
+2\int_{H_m}|f(y)|(1-\chi_\varepsilon(y))\frac{\kappa(m\mid y)}{\kappa(y\mid m)}
\,\mathrm d\mathcal H^{d-1}(y).
\end{aligned}
$$
Letting $\varepsilon\downarrow0$ proves \textup{(ii)}. Taking $f\equiv1$ gives the asserted volume expansion.

Finally, a compact, connected one-dimensional Riemannian
manifold without boundary is isometric to a circle of some
circumference $\ell$.
Let $p$ be the antipodal point of $m$. Then
$$
C_m=H_m=\{p\},
\qquad \kappa(m\mid p)=\kappa(p\mid m)=2.
$$
For $0<\delta<\ell/4$, we have
$\mathcal A_\delta(m)=B_\delta(p)$ and
$\operatorname{vol}(\mathcal A_\delta(m))=2\delta$.
For every $f\in C(\M)$,
$$
\frac1\delta\int_{B_\delta(p)}f\,\mathrm d\operatorname{vol}
\longrightarrow2f(p)
=2\int_{H_m}f(y)
\frac{\kappa(m\mid y)}{\kappa(y\mid m)}
\,\mathrm d\mathcal H^0(y).
$$
Thus both assertions also hold for $d=1$.
\end{proof}

\bibliographystyle{unsrt}
\bibliography{references}

\end{document}